\magnification=1200


\hsize=13.50cm    
\vsize=18cm       
\parindent=12pt   \parskip=0pt     
\pageno=1 

\hoffset=15mm    
\voffset=1cm    
 
\ifnum\mag=\magstep1
\hoffset=-2mm   
\voffset=.8cm   
\fi

\pretolerance=500 \tolerance=1000  \brokenpenalty=5000

\catcode`\@=11

\font\eightrm=cmr8
\font\eighti=cmmi8
\font\eightsy=cmsy8
\font\eightbf=cmbx8
\font\eighttt=cmtt8
\font\eightit=cmti8
\font\eightsl=cmsl8
\font\sevenrm=cmr7
\font\seveni=cmmi7
\font\sevensy=cmsy7
\font\sevenbf=cmbx7

\font\sixrm=cmr6
\font\sixi=cmmi6
\font\sixsy=cmsy6
\font\sixbf=cmbx6

\font\douzebf=cmbx10 at 12pt

\font\tengoth=eufm10
\font\tenbboard=msbm10
\font\eightgoth=eufm10 at 8pt
\font\eightbboard=msbm10 at 8 pt
\font\sevengoth=eufm7
\font\sevenbboard=msbm7
\font\sixgoth=eufm7 at 6 pt
\font\fivegoth=eufm5

\font\tencyr=wncyr10

\font\eightcyr=wncyr10 at 8 pt

\font\sevencyr=wncyr10 at 7 pt

\font\sixcyr=wncyr10 at 6 pt

\skewchar\eighti='177 \skewchar\sixi='177
\skewchar\eightsy='60 \skewchar\sixsy='60

\newfam\gothfam
\newfam\bboardfam
\newfam\cyrfam

\def\tenpoint{%
  \textfont0=\tenrm \scriptfont0=\sevenrm
  \scriptscriptfont0=\fiverm
  \def\rm{\fam\z@\tenrm}%
  \textfont1=\teni  \scriptfont1=\seveni
  \scriptscriptfont1=\fivei
  \def\oldstyle{\fam\@ne\teni}\let\old=\oldstyle
  \textfont2=\tensy \scriptfont2=\sevensy
  \scriptscriptfont2=\fivesy
  \textfont\gothfam=\tengoth
  \scriptfont\gothfam=\sevengoth
  \scriptscriptfont\gothfam=\fivegoth
  \def\goth{\fam\gothfam\tengoth}%
  \textfont\bboardfam=\tenbboard
  \scriptfont\bboardfam=\sevenbboard
  \scriptscriptfont\bboardfam=\sevenbboard
  \def\bb{\fam\bboardfam\tenbboard}%
  \textfont\cyrfam=\tencyr
  \scriptfont\cyrfam=\sevencyr
  \scriptscriptfont\cyrfam=\sixcyr
  \def\cyr{\fam\cyrfam\tencyr}%
  \textfont\itfam=\tenit
  \def\it{\fam\itfam\tenit}%
  \textfont\slfam=\tensl
  \def\sl{\fam\slfam\tensl}%
  \textfont\bffam=\tenbf
  \scriptfont\bffam=\sevenbf
  \scriptscriptfont\bffam=\fivebf
  \def\bf{\fam\bffam\tenbf}%
  \textfont\ttfam=\tentt
  \def\tt{\fam\ttfam\tentt}%
  \abovedisplayskip=12pt plus 3pt minus 9pt
  \belowdisplayskip=\abovedisplayskip
  \abovedisplayshortskip=0pt plus 3pt
  \belowdisplayshortskip=4pt plus 3pt 
  \smallskipamount=3pt plus 1pt minus 1pt
  \medskipamount=6pt plus 2pt minus 2pt
  \bigskipamount=12pt plus 4pt minus 4pt
  \normalbaselineskip=12pt
  \setbox\strutbox=\hbox{\vrule height8.5pt depth3.5pt width0pt}%
  \let\bigf@nt=\tenrm
  \let\smallf@nt=\sevenrm
  \normalbaselines\rm}
  
\def\eightpoint{%
  \textfont0=\eightrm \scriptfont0=\sixrm
  \scriptscriptfont0=\fiverm
  \def\rm{\fam\z@\eightrm}%
  \textfont1=\eighti  \scriptfont1=\sixi
  \scriptscriptfont1=\fivei
  \def\oldstyle{\fam\@ne\eighti}\let\old=\oldstyle
  \textfont2=\eightsy \scriptfont2=\sixsy
  \scriptscriptfont2=\fivesy
  \textfont\gothfam=\eightgoth
  \scriptfont\gothfam=\sixgoth
  \scriptscriptfont\gothfam=\fivegoth
  \def\goth{\fam\gothfam\eightgoth}%
  \textfont\cyrfam=\eightcyr
  \scriptfont\cyrfam=\sixcyr
  \scriptscriptfont\cyrfam=\sixcyr
  \def\cyr{\fam\cyrfam\eightcyr}%
  \textfont\bboardfam=\eightbboard
  \scriptfont\bboardfam=\sevenbboard
  \scriptscriptfont\bboardfam=\sevenbboard
  \def\bb{\fam\bboardfam}%
  \textfont\itfam=\eightit
  \def\it{\fam\itfam\eightit}%
  \textfont\slfam=\eightsl
  \def\sl{\fam\slfam\eightsl}%
  \textfont\bffam=\eightbf
  \scriptfont\bffam=\sixbf
  \scriptscriptfont\bffam=\fivebf
  \def\bf{\fam\bffam\eightbf}%
  \textfont\ttfam=\eighttt
  \def\tt{\fam\ttfam\eighttt}%
  \abovedisplayskip=9pt plus 3pt minus 9pt
  \belowdisplayskip=\abovedisplayskip
  \abovedisplayshortskip=0pt plus 3pt
  \belowdisplayshortskip=3pt plus 3pt 
  \smallskipamount=2pt plus 1pt minus 1pt
  \medskipamount=4pt plus 2pt minus 1pt
  \bigskipamount=9pt plus 3pt minus 3pt
  \normalbaselineskip=9pt
  \setbox\strutbox=\hbox{\vrule height7pt depth2pt width0pt}%
  \let\bigf@nt=\eightrm
  \let\smallf@nt=\sixrm
  \normalbaselines\rm}
\tenpoint


\def\pc#1{\bigf@nt#1\smallf@nt}
\def\pd#1 {{\pc#1} }

\catcode`\;=\active
\def;{\relax\ifhmode\ifdim\lastskip>\z@\unskip\fi
\kern\fontdimen2  -1.2 \fontdimen3 \string;}

\catcode`\:=\active
\def:{\relax\ifhmode\ifdim\lastskip>\z@\unskip\fi\penalty\@M\ 
\fi\string:}

\catcode`\!=\active
\def!{\relax\ifhmode\ifdim\lastskip>\z@
\unskip\fi\kern\fontdimen2  -1.1 \fontdimen3
\string!}

\catcode`\?=\active
\def?{\relax\ifhmode\ifdim\lastskip>\z@
\unskip\fi\kern\fontdimen2  -1.1 \fontdimen3
\string?}

\def\^#1{\if#1i{\accent"5E\i}\else{\accent"5E #1}\fi}
\def\"#1{\if#1i{\accent"7F\i}\else{\accent"7F #1}\fi}

\frenchspacing


\newtoks\auteurcourant
\auteurcourant={\hfil}
\newtoks\titrecourant
\titrecourant={\hfil}

\newtoks\hautpagetitre
\hautpagetitre={\hfil}
\newtoks\baspagetitre
\baspagetitre={\hfil}

\newtoks\hautpagegauche     
\hautpagegauche={\eightpoint\rlap{\folio}\hfil\the
\auteurcourant\hfil}
\newtoks\hautpagedroite     
\hautpagedroite={\eightpoint\hfil\the
\titrecourant\hfil\llap{\folio}}

\newtoks\baspagegauche
\baspagegauche={\hfil} 
\newtoks\baspagedroite
\baspagedroite={\hfil}

\newif\ifpagetitre
\pagetitretrue  


\headline={\ifpagetitre\the\hautpagetitre
\else\ifodd\pageno\the\hautpagedroite\else\the
\hautpagegauche\fi\fi}

\footline={\ifpagetitre\the\baspagetitre\else
\ifodd\pageno\the\baspagedroite\else\the
\baspagegauche\fi\fi
\global\pagetitrefalse}


\def\raggedbottom{\topskip 10pt plus 36pt\r@ggedbottomtrue}


\def\pointir{\unskip . --- \ignorespaces}

\def\Bigbreak{\vskip-\lastskip\bigbreak}
\def\Medbreak{\vskip-\lastskip\medbreak}

\def\ctexte#1\endctexte{%
  \hbox{$\vcenter{\halign{\hfill##\hfill\crcr#1\crcr}}$}}

\long\def\ctitre#1\endctitre{%
    \ifdim\lastskip<24pt\vskip-\lastskip\bigbreak\bigbreak\fi
  		\vbox{\parindent=0pt\leftskip=0pt plus 1fill
          \rightskip=\leftskip
          \parfillskip=0pt\bf#1\par}
    \bigskip\nobreak}

\long\def\section#1\endsection{%
\vskip 0pt plus 3\normalbaselineskip
\penalty-250
\vskip 0pt plus -3\normalbaselineskip
\Bigbreak
\message{[section \string: #1]}{\bf#1\unskip}\pointir}

\long\def\sectiona#1\endsection{%
\vskip 0pt plus 3\normalbaselineskip
\penalty-250
\vskip 0pt plus -3\normalbaselineskip
\Bigbreak
\message{[sectiona \string: #1]}%
{\bf#1}\medskip\nobreak}

\long\def\subsection#1\endsubsection{%
\Medbreak
{\it#1\unskip}\pointir}

\long\def\subsectiona#1\endsubsection{%
\Medbreak
{\it#1}\par\nobreak}

\def\rem#1\endrem{%
\Medbreak
{\it#1\unskip} : }

\def\remp#1\endrem{%
\Medbreak
{\pc #1\unskip}\pointir}

\def\rema#1\endrem{%
\Medbreak
{\it #1}\par\nobreak}

\def\newparwithcolon#1\endnewparwithcolon{
\Medbreak
{#1\unskip} : }

\def\newparwithpointir#1\endnewparwithpointir{
\Medbreak
{#1\unskip}\pointir}

\def\newpara#1\endnewpar{
\Medbreak
{#1\unskip}\smallskip\nobreak}


\long\def\th#1 #2\enonce#3\endth{%
   \Medbreak
   {\pc#1} {#2\unskip}\pointir{\it #3}\medskip}

\long\def\tha#1 #2\enonce#3\endth{%
   \Medbreak
   {\pc#1} {#2\unskip}\par\nobreak{\it #3}\medskip}

\long\def\Th#1 #2 #3\enonce#4\endth{%
   \Medbreak
   #1 {\pc#2} {#3\unskip}\pointir{\it #4}\medskip}

\long\def\Tha#1 #2 #3\enonce#4\endth{%
   \Medbreak
   #1 {\pc#2} #3\par\nobreak{\it #4}\medskip}

\def\decale#1{\smallbreak\hskip 28pt\llap{#1}\kern 5pt}
\def\decaledecale#1{\smallbreak\hskip 34pt\llap{#1}\kern 5pt}
\def\puce{\smallbreak\hskip 6pt{$\scriptstyle\bullet$}\kern 5pt}


\def\displaylinesno#1{\displ@y\halign{
\hbox to\displaywidth{$\@lign\hfil\displaystyle##\hfil$}&
\llap{$##$}\crcr#1\crcr}}

\def\ldisplaylinesno#1{\displ@y\halign{ 
\hbox to\displaywidth{$\@lign\hfil\displaystyle##\hfil$}&
\kern-\displaywidth\rlap{$##$}\tabskip\displaywidth\crcr#1\crcr}}

\def\eqalign#1{\null\,\vcenter{\openup\jot\m@th\ialign{
\strut\hfil$\displaystyle{##}$&$\displaystyle{{}##}$\hfil
&&\quad\strut\hfil$\displaystyle{##}$&$\displaystyle{{}##}$\hfil
\crcr#1\crcr}}\,}


\def\system#1{\left\{\null\,\vcenter{\openup1\jot\m@th
\ialign{\strut$##$&\hfil$##$&$##$\hfil&&
        \enskip$##$\enskip&\hfil$##$&$##$\hfil\crcr#1\crcr}}
        \right.}

\let\@ldmessage=\message

\def\message#1{{\def\pc{\string\pc\space}%
                \def\'{\string'}\def\`{\string`}%
                \def\^{\string^}\def\"{\string"}%
                \@ldmessage{#1}}}

\def\({{\rm (}}
\def\){{\rm )}}

\def\up#1{\raise 1ex\hbox{\smallf@nt#1}}


\def\qed{\raise -2pt\hbox{\vrule\vbox to 10pt{\hrule width 4pt
                 \vfill\hrule}\vrule}}

\def\cqfd{\unskip\penalty 500\quad\qed\medbreak}

\def\virg{\raise .4ex\hbox{,}}   


\def\build#1_#2^#3{\mathrel{
\mathop{\kern 0pt#1}\limits_{#2}^{#3}}}

\def\boxit#1#2{%
\setbox1=\hbox{\kern#1{#2}\kern#1}%
\dimen1=\ht1 \advance\dimen1 by #1 \dimen2=\dp1
\advance\dimen2 by #1 
\setbox1=\hbox{\vrule height\dimen1
depth\dimen2\box1\vrule}%
\setbox1=\vbox{\hrule\box1\hrule}%
\advance\dimen1 by .6pt \ht1=\dimen1 
\advance\dimen2 by .6pt \dp1=\dimen2  \box1\relax}

\catcode`\@=12

\showboxbreadth=-1  \showboxdepth=-1


\def\Grille{\setbox13=\vbox to 5\unitlength{\hrule width 109mm
\vfill} 
\setbox13=\vbox to 65mm
{\offinterlineskip\leaders\copy13\vfill\kern-1pt\hrule} 
\setbox14=\hbox to 5\unitlength{\vrule height 65mm\hfill} 
\setbox14=\hbox to 109mm{\leaders\copy14\hfill\kern-2mm
\vrule height 65mm}
\ht14=0pt\dp14=0pt\wd14=0pt \setbox13=\vbox to 0pt
{\vss\box13\offinterlineskip\box14} \wd13=0pt\box13}

\def\arrow(#1,#2)\dir(#3,#4)\length#5{%
\noalign{\leftput(#1,#2){\vector(#3,#4){#5}}}}

\def\ligne(#1,#2)\dir(#3,#4)\length#5{%
\noalign{\leftput(#1,#2){\lline(#3,#4){#5}}}}

\def\put(#1,#2)#3{\noalign{\setbox1=\hbox{%
    \kern #1\unitlength
    \raise #2\unitlength\hbox{$#3$}}%
    \ht1=0pt \wd1=0pt \dp1=0pt\box1}}

\def\diagram#1{\def\normalbaselines{\baselineskip=0pt
\lineskip=5pt}\matrix{#1}}

\def\maprightover#1{\smash{\mathop{\longrightarrow}
\limits^{#1}}}

\def\maprightunder#1{\smash{\mathop{\longrightarrow}
\limits_{#1}}}

\def\mapdownleft#1{\llap{$\vcenter
{\hbox{$\scriptstyle#1$}}$}
\Big\downarrow}

\def\mapdownright#1{\Big\downarrow
\rlap{$\vcenter{\hbox{$\scriptstyle#1$}}$}}

\def\longmaprightover#1#2{\smash{\mathop{\hbox to#2
{\rightarrowfill}}\limits^{\scriptstyle#1}}}

\def\longmapleftover#1#2{\smash{\mathop{\hbox to#2
{\leftarrowfill}}\limits^{\scriptstyle#1}}}

\def\longmaprightunder#1#2{\smash{\mathop{\hbox to#2
{\rightarrowfill}}\limits_{\scriptstyle#1}}}

\def\longmaplefttunder#1#2{\smash{\mathop{\hbox to#2
{\leftarrowfill}}\limits_{\scriptstyle#1}}}

\def\longhookrightarrowover#1#2{\smash{\mathop{\lhook\joinrel
\mathrel{\hbox to #2{\rightarrowfill}}}\limits^{\scriptstyle#1}}}

\def\longhookrightarrowunder#1#2{\smash{\mathop{\lhook\joinrel
\mathrel{\hbox to #2{\rightarrowfill}}}\limits_{\scriptstyle#1}}}

\def\longhookleftarrowover#1#2{\smash{\mathop{{\hbox to #2
{\leftarrowfill}}\joinrel\kern -0.9mm\mathrel\rhook}
\limits^{\scriptstyle#1}}}

\def\longhookleftarrowunder#1#2{\smash{\mathop{{\hbox to #2
{\leftarrowfill}}\joinrel\kern -0.9mm\mathrel\rhook}
\limits_{\scriptstyle#1}}}

\def\longtwoheadrightarrowover#1#2{\smash{\mathop{{\hbox to #2
{\rightarrowfill}}\kern -3.25mm\joinrel\mathrel\rightarrow}
\limits^{\scriptstyle#1}}}

\def\longtwoheadrightarrowunder#1#2{\smash{\mathop{{\hbox to #2
{\rightarrowfill}}\kern -3.25mm\joinrel\mathrel\rightarrow}
\limits_{\scriptstyle#1}}}

\def\longtwoheadleftarrowover#1#2{\smash{\mathop{\joinrel\mathrel
\leftarrow\kern -3.8mm{\hbox to #2{\leftarrowfill}}}
\limits^{\scriptstyle#1}}}

\def\longtwoheadleftarrowunder#1#2{\smash{\mathop{\joinrel\mathrel
\leftarrow\kern -3.8mm{\hbox to #2{\leftarrowfill}}}
\limits_{\scriptstyle#1}}}

\def\isomorphism{\buildrel\sim\over\longrightarrow }

\def\longmapsto#1{\mapstochar\mathrel{\joinrel
\kern-0.2mm\hbox to #1mm{\rightarrowfill}}}

\message{`lline' & `vector' macros from LaTeX}
\catcode`@=11
\def\{{\relax\ifmmode\lbrace\else$\lbrace$\fi}
\def\}{\relax\ifmmode\rbrace\else$\rbrace$\fi}
\def\newcount{\alloc@0\count\countdef\insc@unt}
\def\newdimen{\alloc@1\dimen\dimendef\insc@unt}
\def\newwrite{\alloc@7\write\chardef\sixt@@n}

\newwrite\@unused
\newcount\@tempcnta
\newcount\@tempcntb
\newdimen\@tempdima
\newdimen\@tempdimb
\newbox\@tempboxa

\def\@spaces{\space\space\space\space}
\def\@whilenoop#1{}
\def\@whiledim#1\do #2{\ifdim #1\relax#2\@iwhiledim{#1\relax#2}\fi}
\def\@iwhiledim#1{\ifdim #1\let\@nextwhile=\@iwhiledim
        \else\let\@nextwhile=\@whilenoop\fi\@nextwhile{#1}}
\def\@badlinearg{\@latexerr{Bad \string\line\space or \string\vector
   \space argument}}
\def\@latexerr#1#2{\begingroup
\edef\@tempc{#2}\expandafter\errhelp\expandafter{\@tempc}%

\def\@eha{Your command was ignored.
^^JType \space I <command> <return> \space to replace it
  with another command,^^Jor \space <return> \space to continue without
it.} 
\def\@ehb{You've lost some text. \space \@ehc}
\def\@ehc{Try typing \space <return>
  \space to proceed.^^JIf that doesn't work, type \space X <return>
  \space to quit.}
\def\@ehd{You're in trouble here.  \space\@ehc}

\typeout{LaTeX error. \space See LaTeX manual for explanation.^^J
 \space\@spaces\@spaces\@spaces Type \space H <return> \space for
 immediate help.}\errmessage{#1}\endgroup}
\def\typeout#1{{\let\protect\string\immediate\write\@unused{#1}}}

\font\tenln    = line10
\font\tenlnw   = linew10

\newdimen\@wholewidth
\newdimen\@halfwidth
\newdimen\unitlength 

\unitlength =1pt

\def\thinlines{\let\@linefnt\tenln \let\@circlefnt\tencirc
  \@wholewidth\fontdimen8\tenln \@halfwidth .5\@wholewidth}
\def\thicklines{\let\@linefnt\tenlnw \let\@circlefnt\tencircw
  \@wholewidth\fontdimen8\tenlnw \@halfwidth .5\@wholewidth}

\def\linethickness#1{\@wholewidth #1\relax \@halfwidth .5
\@wholewidth}

\newif\if@negarg

\def\lline(#1,#2)#3{\@xarg #1\relax \@yarg #2\relax
\@linelen=#3\unitlength
\ifnum\@xarg =0 \@vline
  \else \ifnum\@yarg =0 \@hline \else \@sline\fi
\fi}

\def\@sline{\ifnum\@xarg< 0 \@negargtrue \@xarg -\@xarg
  \@yyarg -\@yarg
  \else \@negargfalse \@yyarg \@yarg \fi
\ifnum \@yyarg >0 \@tempcnta\@yyarg \else \@tempcnta -
  \@yyarg \fi
\ifnum\@tempcnta>6 \@badlinearg\@tempcnta0 \fi
\setbox\@linechar\hbox{\@linefnt\@getlinechar(\@xarg,\@yyarg)}%
\ifnum \@yarg >0 \let\@upordown\raise \@clnht\z@
   \else\let\@upordown\lower \@clnht \ht\@linechar\fi
\@clnwd=\wd\@linechar
\if@negarg \hskip -\wd\@linechar \def\@tempa{\hskip -2\wd
  \@linechar}
  \else \let\@tempa\relax \fi
\@whiledim \@clnwd <\@linelen \do
  {\@upordown\@clnht\copy\@linechar
   \@tempa
   \advance\@clnht \ht\@linechar
   \advance\@clnwd \wd\@linechar}%
\advance\@clnht -\ht\@linechar
\advance\@clnwd -\wd\@linechar
\@tempdima\@linelen\advance\@tempdima -\@clnwd
\@tempdimb\@tempdima\advance\@tempdimb -\wd\@linechar
\if@negarg \hskip -\@tempdimb \else \hskip \@tempdimb \fi
\multiply\@tempdima \@m\@tempcnta \@tempdima \@tempdima
\wd\@linechar \divide\@tempcnta \@tempdima
\@tempdima \ht\@linechar \multiply\@tempdima \@tempcnta
\divide\@tempdima \@m
\advance\@clnht \@tempdima
\ifdim \@linelen <\wd\@linechar
   \hskip \wd\@linechar
  \else\@upordown\@clnht\copy\@linechar\fi}

\def\@hline{\ifnum \@xarg <0 \hskip -\@linelen \fi
\vrule height \@halfwidth depth \@halfwidth width \@linelen
\ifnum \@xarg <0 \hskip -\@linelen \fi}

\def\@getlinechar(#1,#2){\@tempcnta#1\relax
\multiply\@tempcnta 8\advance\@tempcnta -9
\ifnum #2>0 \advance\@tempcnta #2\relax
 \else\advance\@tempcnta -#2\relax\advance\@tempcnta 64 \fi
\char\@tempcnta}

\def\vector(#1,#2)#3{\@xarg #1\relax \@yarg #2\relax
\@linelen=#3\unitlength
\ifnum\@xarg =0 \@vvector
  \else \ifnum\@yarg =0 \@hvector \else \@svector\fi
\fi}

\def\@hvector{\@hline\hbox to 0pt{\@linefnt
\ifnum \@xarg <0 \@getlarrow(1,0)\hss\else
    \hss\@getrarrow(1,0)\fi}}

\def\@vvector{\ifnum \@yarg <0 \@downvector \else \@upvector \fi}

\def\@svector{\@sline\@tempcnta\@yarg
\ifnum\@tempcnta <0 \@tempcnta=-\@tempcnta\fi
\ifnum\@tempcnta <5
  \hskip -\wd\@linechar
  \@upordown\@clnht \hbox{\@linefnt  \if@negarg
  \@getlarrow(\@xarg,\@yyarg)
  \else \@getrarrow(\@xarg,\@yyarg) \fi}%
\else\@badlinearg\fi}

\def\@getlarrow(#1,#2){\ifnum #2 =\z@ \@tempcnta='33\else
\@tempcnta=#1\relax\multiply\@tempcnta \sixt@@n
\advance\@tempcnta -9 \@tempcntb=#2\relax
\multiply\@tempcntb \tw@
\ifnum \@tempcntb >0 \advance\@tempcnta \@tempcntb\relax
 \else\advance\@tempcnta -\@tempcntb\advance\@tempcnta 64
\fi\fi
\char\@tempcnta}

\def\@getrarrow(#1,#2){\@tempcntb=#2\relax
\ifnum\@tempcntb < 0 \@tempcntb=-\@tempcntb\relax\fi
\ifcase \@tempcntb\relax \@tempcnta='55 \or
\ifnum #1<3 \@tempcnta=#1\relax\multiply\@tempcnta
24 \advance\@tempcnta -6 \else \ifnum #1=3 \@tempcnta=49
\else\@tempcnta=58 \fi\fi\or
\ifnum #1<3 \@tempcnta=#1\relax\multiply\@tempcnta
24 \advance\@tempcnta -3 \else \@tempcnta=51\fi\or
\@tempcnta=#1\relax\multiply\@tempcnta
\sixt@@n \advance\@tempcnta -\tw@ \else
\@tempcnta=#1\relax\multiply\@tempcnta
\sixt@@n \advance\@tempcnta 7 \fi
\ifnum #2<0 \advance\@tempcnta 64 \fi
\char\@tempcnta}

\def\@vline{\ifnum \@yarg <0 \@downline \else \@upline\fi}

\def\@upline{\hbox to \z@{\hskip -\@halfwidth \vrule
  width \@wholewidth height \@linelen depth \z@\hss}}

\def\@downline{\hbox to \z@{\hskip -\@halfwidth \vrule
  width \@wholewidth height \z@ depth \@linelen \hss}}

\def\@upvector{\@upline\setbox\@tempboxa
     \hbox{\@linefnt\char'66}\raise
     \@linelen \hbox to\z@{\lower \ht\@tempboxa
     \box\@tempboxa\hss}}

\def\@downvector{\@downline\lower \@linelen
      \hbox to \z@{\@linefnt\char'77\hss}}

\thinlines

\newcount\@xarg
\newcount\@yarg
\newcount\@yyarg
\newcount\@multicnt
\newdimen\@xdim
\newdimen\@ydim
\newbox\@linechar
\newdimen\@linelen
\newdimen\@clnwd
\newdimen\@clnht
\newdimen\@dashdim
\newbox\@dashbox
\newcount\@dashcnt
 \catcode`@=12

\newbox\tbox
\newbox\tboxa

\def\leftzer#1{\setbox\tbox=\hbox to 0pt{#1\hss}%
     \ht\tbox=0pt \dp\tbox=0pt \box\tbox}

\def\rightzer#1{\setbox\tbox=\hbox to 0pt{\hss #1}%
     \ht\tbox=0pt \dp\tbox=0pt \box\tbox}

\def\centerzer#1{\setbox\tbox=\hbox to 0pt{\hss #1\hss}%
     \ht\tbox=0pt \dp\tbox=0pt \box\tbox}

%
\def\image(#1,#2)#3{\vbox to #1{\offinterlineskip
    \vss #3 \vskip #2}}

\def\leftput(#1,#2)#3{\setbox\tboxa=\hbox{%
    \kern #1\unitlength
    \raise #2\unitlength\hbox{\leftzer{#3}}}%
    \ht\tboxa=0pt \wd\tboxa=0pt \dp\tboxa=0pt\box\tboxa}

\def\rightput(#1,#2)#3{\setbox\tboxa=\hbox{%
    \kern #1\unitlength
    \raise #2\unitlength\hbox{\rightzer{#3}}}%
    \ht\tboxa=0pt \wd\tboxa=0pt \dp\tboxa=0pt\box\tboxa}

\def\centerput(#1,#2)#3{\setbox\tboxa=\hbox{%
    \kern #1\unitlength
    \raise #2\unitlength\hbox{\centerzer{#3}}}%
    \ht\tboxa=0pt \wd\tboxa=0pt \dp\tboxa=0pt\box\tboxa}

\unitlength=1mm

\expandafter\ifx\csname amssym.def\endcsname\relax \else
\endinput\fi
%
\expandafter\edef\csname amssym.def\endcsname{%
       \catcode`\noexpand\@=\the\catcode`\@\space}
\catcode`\@=11
%

\def\undefine#1{\let#1\undefined}
\def\newsymbol#1#2#3#4#5{\let\next@\relax
 \ifnum#2=\@ne\let\next@\msafam@\else
 \ifnum#2=\tw@\let\next@\msbfam@\fi\fi
 \mathchardef#1="#3\next@#4#5}
\def\mathhexbox@#1#2#3{\relax
 \ifmmode\mathpalette{}{\m@th\mathchar"#1#2#3}%
 \else\leavevmode\hbox{$\m@th\mathchar"#1#2#3$}\fi}
\def\hexnumber@#1{\ifcase#1 0\or 1\or 2\or 3
\or 4\or 5\or 6\or 7\or 8\or
 9\or A\or B\or C\or D\or E\or F\fi}

\font\tenmsa=msam10
\font\sevenmsa=msam7
\font\fivemsa=msam5
\newfam\msafam
\textfont\msafam=\tenmsa
\scriptfont\msafam=\sevenmsa
\scriptscriptfont\msafam=\fivemsa
\edef\msafam@{\hexnumber@\msafam}
\mathchardef\dabar@"0\msafam@39
\def\dashrightarrow{\mathrel{\dabar@\dabar@\mathchar"0
\msafam@4B}}
\def\dashleftarrow{\mathrel{\mathchar"0\msafam@4C
\dabar@\dabar@}}

\def\ulcorner{\delimiter"4\msafam@70\msafam@70 }
\def\urcorner{\delimiter"5\msafam@71\msafam@71 }
\def\llcorner{\delimiter"4\msafam@78\msafam@78 }
\def\lrcorner{\delimiter"5\msafam@79\msafam@79 }
\def\yen{{\mathhexbox@\msafam@55}}
\def\checkmark{{\mathhexbox@\msafam@58}}
\def\circledR{{\mathhexbox@\msafam@72}}
\def\maltese{{\mathhexbox@\msafam@7A}}

\font\tenmsb=msbm10
\font\sevenmsb=msbm7
\font\fivemsb=msbm5
\newfam\msbfam
\textfont\msbfam=\tenmsb
\scriptfont\msbfam=\sevenmsb
\scriptscriptfont\msbfam=\fivemsb
\edef\msbfam@{\hexnumber@\msbfam}
\def\Bbb#1{{\fam\msbfam\relax#1}}
\def\widehat#1{\setbox\z@\hbox{$\m@th#1$}%
 \ifdim\wd\z@>\tw@ em\mathaccent"0\msbfam@5B{#1}%
 \else\mathaccent"0362{#1}\fi}
\def\widetilde#1{\setbox\z@\hbox{$\m@th#1$}%
 \ifdim\wd\z@>\tw@ em\mathaccent"0\msbfam@5D{#1}%
 \else\mathaccent"0365{#1}\fi}
\font\teneufm=eufm10
\font\seveneufm=eufm7
\font\fiveeufm=eufm5
\newfam\eufmfam
\textfont\eufmfam=\teneufm
\scriptfont\eufmfam=\seveneufm
\scriptscriptfont\eufmfam=\fiveeufm
\def\frak#1{{\fam\eufmfam\relax#1}}
\let\goth\frak

\csname amssym.def\endcsname

\expandafter\ifx\csname pre amssym.tex at\endcsname\relax \else 
\endinput\fi
\expandafter\chardef\csname pre amssym.tex at\endcsname=\the
\catcode`\@
\catcode`\@=11
\begingroup\ifx\undefined\newsymbol \else\def\input#1
{\endgroup}\fi
\input amssym.def \relax
\newsymbol\boxdot 1200
\newsymbol\boxplus 1201
\newsymbol\boxtimes 1202
\newsymbol\square 1003
\newsymbol\blacksquare 1004
\newsymbol\centerdot 1205
\newsymbol\lozenge 1006
\newsymbol\blacklozenge 1007
\newsymbol\circlearrowright 1308
\newsymbol\circlearrowleft 1309
\undefine\rightleftharpoons
\newsymbol\rightleftharpoons 130A
\newsymbol\leftrightharpoons 130B
\newsymbol\boxminus 120C
\newsymbol\Vdash 130D
\newsymbol\Vvdash 130E
\newsymbol\vDash 130F
\newsymbol\twoheadrightarrow 1310
\newsymbol\twoheadleftarrow 1311
\newsymbol\leftleftarrows 1312
\newsymbol\rightrightarrows 1313
\newsymbol\upuparrows 1314
\newsymbol\downdownarrows 1315
\newsymbol\upharpoonright 1316
 
\newsymbol\downharpoonright 1317
\newsymbol\upharpoonleft 1318
\newsymbol\downharpoonleft 1319
\newsymbol\rightarrowtail 131A
\newsymbol\leftarrowtail 131B
\newsymbol\leftrightarrows 131C
\newsymbol\rightleftarrows 131D
\newsymbol\Lsh 131E
\newsymbol\Rsh 131F
\newsymbol\rightsquigarrow 1320
\newsymbol\leftrightsquigarrow 1321
\newsymbol\looparrowleft 1322
\newsymbol\looparrowright 1323
\newsymbol\circeq 1324
\newsymbol\succsim 1325
\newsymbol\gtrsim 1326
\newsymbol\gtrapprox 1327
\newsymbol\multimap 1328
\newsymbol\therefore 1329
\newsymbol\because 132A
\newsymbol\doteqdot 132B
 
\newsymbol\triangleq 132C
\newsymbol\precsim 132D
\newsymbol\lesssim 132E
\newsymbol\lessapprox 132F
\newsymbol\eqslantless 1330
\newsymbol\eqslantgtr 1331
\newsymbol\curlyeqprec 1332
\newsymbol\curlyeqsucc 1333
\newsymbol\preccurlyeq 1334
\newsymbol\leqq 1335
\newsymbol\leqslant 1336
\newsymbol\lessgtr 1337
\newsymbol\backprime 1038
\newsymbol\risingdotseq 133A
\newsymbol\fallingdotseq 133B
\newsymbol\succcurlyeq 133C
\newsymbol\geqq 133D
\newsymbol\geqslant 133E
\newsymbol\gtrless 133F
\newsymbol\sqsubset 1340
\newsymbol\sqsupset 1341
\newsymbol\vartriangleright 1342
\newsymbol\vartriangleleft 1343
\newsymbol\trianglerighteq 1344
\newsymbol\trianglelefteq 1345
\newsymbol\bigstar 1046
\newsymbol\between 1347
\newsymbol\blacktriangledown 1048
\newsymbol\blacktriangleright 1349
\newsymbol\blacktriangleleft 134A
\newsymbol\vartriangle 134D
\newsymbol\blacktriangle 104E
\newsymbol\triangledown 104F
\newsymbol\eqcirc 1350
\newsymbol\lesseqgtr 1351
\newsymbol\gtreqless 1352
\newsymbol\lesseqqgtr 1353
\newsymbol\gtreqqless 1354
\newsymbol\Rrightarrow 1356
\newsymbol\Lleftarrow 1357
\newsymbol\veebar 1259
\newsymbol\barwedge 125A
\newsymbol\doublebarwedge 125B
\undefine\angle
\newsymbol\angle 105C
\newsymbol\measuredangle 105D
\newsymbol\sphericalangle 105E
\newsymbol\varpropto 135F
\newsymbol\smallsmile 1360
\newsymbol\smallfrown 1361
\newsymbol\Subset 1362
\newsymbol\Supset 1363
\newsymbol\Cup 1264
 
\newsymbol\Cap 1265
 
\newsymbol\curlywedge 1266
\newsymbol\curlyvee 1267
\newsymbol\leftthreetimes 1268
\newsymbol\rightthreetimes 1269
\newsymbol\subseteqq 136A
\newsymbol\supseteqq 136B
\newsymbol\bumpeq 136C
\newsymbol\Bumpeq 136D
\newsymbol\lll 136E
 
\newsymbol\ggg 136F
 
\newsymbol\circledS 1073
\newsymbol\pitchfork 1374
\newsymbol\dotplus 1275
\newsymbol\backsim 1376
\newsymbol\backsimeq 1377
\newsymbol\complement 107B
\newsymbol\intercal 127C
\newsymbol\circledcirc 127D
\newsymbol\circledast 127E
\newsymbol\circleddash 127F
\newsymbol\lvertneqq 2300
\newsymbol\gvertneqq 2301
\newsymbol\nleq 2302
\newsymbol\ngeq 2303
\newsymbol\nless 2304
\newsymbol\ngtr 2305
\newsymbol\nprec 2306
\newsymbol\nsucc 2307
\newsymbol\lneqq 2308
\newsymbol\gneqq 2309
\newsymbol\nleqslant 230A
\newsymbol\ngeqslant 230B
\newsymbol\lneq 230C
\newsymbol\gneq 230D
\newsymbol\npreceq 230E
\newsymbol\nsucceq 230F
\newsymbol\precnsim 2310
\newsymbol\succnsim 2311
\newsymbol\lnsim 2312
\newsymbol\gnsim 2313
\newsymbol\nleqq 2314
\newsymbol\ngeqq 2315
\newsymbol\precneqq 2316
\newsymbol\succneqq 2317
\newsymbol\precnapprox 2318
\newsymbol\succnapprox 2319
\newsymbol\lnapprox 231A
\newsymbol\gnapprox 231B
\newsymbol\nsim 231C
\newsymbol\ncong 231D
\newsymbol\diagup 201E
\newsymbol\diagdown 201F
\newsymbol\varsubsetneq 2320
\newsymbol\varsupsetneq 2321
\newsymbol\nsubseteqq 2322
\newsymbol\nsupseteqq 2323
\newsymbol\subsetneqq 2324
\newsymbol\supsetneqq 2325
\newsymbol\varsubsetneqq 2326
\newsymbol\varsupsetneqq 2327
\newsymbol\subsetneq 2328
\newsymbol\supsetneq 2329
\newsymbol\nsubseteq 232A
\newsymbol\nsupseteq 232B
\newsymbol\nparallel 232C
\newsymbol\nmid 232D
\newsymbol\nshortmid 232E
\newsymbol\nshortparallel 232F
\newsymbol\nvdash 2330
\newsymbol\nVdash 2331
\newsymbol\nvDash 2332
\newsymbol\nVDash 2333
\newsymbol\ntrianglerighteq 2334
\newsymbol\ntrianglelefteq 2335
\newsymbol\ntriangleleft 2336
\newsymbol\ntriangleright 2337
\newsymbol\nleftarrow 2338
\newsymbol\nrightarrow 2339
\newsymbol\nLeftarrow 233A
\newsymbol\nRightarrow 233B
\newsymbol\nLeftrightarrow 233C
\newsymbol\nleftrightarrow 233D
\newsymbol\divideontimes 223E
\newsymbol\varnothing 203F
\newsymbol\nexists 2040
\newsymbol\Finv 2060
\newsymbol\Game 2061
\newsymbol\mho 2066
\newsymbol\eth 2067
\newsymbol\eqsim 2368
\newsymbol\beth 2069
\newsymbol\gimel 206A
\newsymbol\daleth 206B
\newsymbol\lessdot 236C
\newsymbol\gtrdot 236D
\newsymbol\ltimes 226E
\newsymbol\rtimes 226F
\newsymbol\shortmid 2370
\newsymbol\shortparallel 2371
\newsymbol\smallsetminus 2272
\newsymbol\thicksim 2373
\newsymbol\thickapprox 2374
\newsymbol\approxeq 2375
\newsymbol\succapprox 2376
\newsymbol\precapprox 2377
\newsymbol\curvearrowleft 2378
\newsymbol\curvearrowright 2379
\newsymbol\digamma 207A
\newsymbol\varkappa 207B
\newsymbol\Bbbk 207C
\newsymbol\hslash 207D
\undefine\hbar
\newsymbol\hbar 207E
\newsymbol\backepsilon 237F
\catcode`\@=\csname pre amssym.tex at\endcsname

\centerline{\douzebf A geometric approach to the fundamental lemma}
\vskip 1mm
\centerline{\douzebf for unitary groups}
\vskip 10mm
\centerline{G. Laumon and M. Rapoport}
\vskip 20mm

\centerline{\bf 0. Introduction}
\vskip 5mm
Let $n_{1},n_{2}$ be two positive integers, let $G=U(n_{1}+n_{2})$ be 
``the'' unramified unitary group in $n_{1}+n_{2}$ variables over a non 
archimedean local field $F$ and let $H$ be the elliptic endoscopic 
group $U(n_{1})\times U(n_{2})$ for $G$.  Let $K$ be a hyperspecial 
compact open subgroup of $G(F)$ and let $K^H$ be a hyperspecial 
compact open subgroup of $H(F)$.  Then, for any regular semisimple 
element $\gamma\in G(F)$ which comes from an elliptic semisimple 
element $\delta$ in $H(F)$ Langlands and Shelstad ([La-Sh]) have 
defined the $\kappa$-orbital integral $O_\gamma^\kappa (1_{K})$, the 
stable orbital integral $SO_\delta^H(1_{K^H})$ and the transfer factor 
$\Delta (\gamma ,\delta)$ and they have conjectured that
$$
SO_\delta^H (1_{K^H})=\Delta (\delta ,\gamma )O_\gamma^\kappa (1_{K}).  
\leqno (\ast )
$$
This relation is called the fundamental lemma for the pair $(G,H)$.  
It has been proved by Labesse and Langlands in the particular case 
$n_{1}= n_{2}=1$ ([La-La]) and by Kottwitz in the particular case 
$n_{1}=1$ and $n_{2}=2$ ([Ko]).
\vskip 2mm

In this paper we restrict ourselves to the case where $F$ is of equal 
characteristic $p>0$ and we consider the relation $(\ast )$ from a 
geometric point of view.  The restriction to the equal characteristic 
case is more or less equivalent to considering the unequal 
characteristic case under the hypothesis that the residual 
characteristic is ``large enough with respect to the element 
$\delta$''.  In the equal characteristic case the geometric 
interpretation of the orbital integrals is straightforward and does 
not involve Witt vector schemes.  Moreover we restrict ourselves to 
the ``totally ramified case'': we fix totally ramified separable 
extensions $E_{1}$ and $E_{2}$ of $F$ of degree $n_{1}$ and $n_{2}$ 
respectively, we embed the corresponding elliptic torus 
$T=E_{1}'^{1}\times E_{2}'^{1}$ as a maximal torus in both $G$ and in 
$H$ and we take as elements $\gamma$ and $\delta$ the images by these 
embeddings of some element $(\gamma_{1},\gamma_{2})\in T$.  Here 
$E_{i}'$ is the unramified quadratic extension of $E_{i}$ and 
$E_{i}'^{1}\subset E_{i}'^{\times}$ is the subgroup of elements of 
norm $1$ with respect to $E_{i}$.

We use the standard computation of the orbital integrals as numbers of 
selfdual lattices which are fixed by unitary automorphisms in certain 
hermitian vector spaces (see [Ko]).  In this computation the 
$\kappa$-orbital integral appears as the difference between the number 
of selfdual lattices for two different hermitian forms $\Phi^{+}$ and 
$\Phi^{-}$ on $E_{1}'\oplus E_{2}'$ which are fixed by the 
multiplication by $(\gamma_{1},\gamma_{2})$.  The stable orbital 
integral is equal to the product of the number of selfdual lattices in 
$E_{1}'$ and in $E_{2}'$ which are fixed by the multiplication by 
$\gamma_{1}$ and $\gamma_{2}$ respectively.

The transfer factor $\Delta (\delta ,\gamma )$ has been computed by 
Waldspurger for classical groups.  In our case we simply have
$$
\Delta (\delta ,\gamma )=(-1)^{r}q^{-r}
$$
where $q$ is the number of elements in the residue field $k$ and $r$ 
denotes the valuation of the resultant of the minimal polynomials of 
$\gamma_{1}$ and $\gamma_{2}$.
\vskip 2mm

We now explain the contents of this paper.  In a first step (Part I) 
we construct schemes ${\cal X}^{+}$, ${\cal X}^{-}$ and ${\cal 
Y}_{1}$, ${\cal Y}_{2}$ over the residue field $k$ whose $k$-rational 
points are in bijection with the sets of lattices in question.  We thus 
have
$$
O_{\gamma}^{\kappa}(1_{K}) =|{\cal X}^{+}(k)|-|{\cal X}^{-}(k)|
$$
and
$$
SO^{H}_{\delta}(1_{K^{H}}) =|{\cal Y}_{1}(k)|\cdot |{\cal Y}_{2}(k)|.
$$
It turns out that ${\cal Y}_{1}$ and ${\cal Y}_{2}$ are projective 
schemes (closed subschemes of Grassmannians), whereas ${\cal X}^{+}$ 
and ${\cal X}^{-}$ are only locally of finite type.  However ${\cal 
X}^{+}$ and ${\cal X}^{-}$ carry natural actions of ${\Bbb Z}$ such 
that the quotients ${\cal X}^{+}/{\Bbb Z}$ and ${\cal X}^{-}/{\Bbb Z}$ 
are representable by projective schemes over $k$.  Of course, this 
construction may be viewed as a special case of the construction of 
Kazhdan and Lusztig [Ka-Lu].

The closed subschemes ${\cal X}^{+}$ and ${\cal X}^{-}$ contain 
canonical closed subschemes which are projective schemes and contain 
all $k$-rational points.  In order to take into account the sign 
$(-1)^{r}$ of the transfer factor it is convenient to denote these 
subschemes ${\cal X}\subset {\cal X}^{+}$ and ${\cal X}'\subset {\cal 
X}^{-}$ if $r$ is even, and ${\cal X}'\subset {\cal X}^{+}$ and ${\cal 
X}\subset {\cal X}^{-}$ if $r$ is odd.  The geometric version of the 
conjecture of Langlands and Shelstad establishes a close relation 
between the number of points of the schemes introduced above for any 
finite extension of $k$.  Namely, we conjecture that for any extension 
$k_f$ of finite degree $f$ of $k$,
$$
|{\cal X}(k_f)|- |{\cal X}'(k_f)|=q^{fr}\cdot |{\cal Y}_{1}(k_f)|\cdot 
|{\cal Y}_{2}(k_f)|.\leqno{(\ast\ast)}
$$
By the above remarks, the relation $(*)$ is the particular case $f=1$ 
of $(**)$.

The main result of this paper (Part II) is the proof of this 
conjecture for extensions of even degree of $k$, i.e.  extensions of 
the quadratic extension $k'$ of $k$.  For this we construct a 
partition (in fact, two such partitions, interchanged by the Frobenius 
morphism over $k$) of ${\cal X}^{\pm}\otimes_kk'$ into locally closed 
subsets which are vector bundles of rank $r$ over $({\cal Y}_{1}\times 
{\cal Y}_{2}) \otimes_kk'$.  This last assertion is proved by 
interpreting our lattices as coherent modules on germs of singular 
curves contained in ${\rm Spec}(k'[[T_{1}, T_{2}]])$ and using the 
interpretation of $r$ as the intersection multiplicity of these 
curves.  At this point we make use of a fundamental result of Deligne 
on intersection multiplicities ``with weights''.  We then construct a 
closed embedding (in fact, two such embeddings interchanged by the 
Frobenius morphism over $k$) of ${\cal X}'\otimes_kk'$ into ${\cal 
X}\otimes_kk'$ such that the complement is one piece of the above 
partition of ${\cal X}^{\pm}\otimes_kk'$.  Our main result follows now 
by a simple counting argument.

In the final part (Part III) we explain a possible approach to the 
descent from $k'$ to $k$.  Whereas the theory over $k'$ is essentially 
elementary, we envisage the use of $\ell$-adic cohomology for this 
descent problem.  Briefly put, even though the vector bundle structure 
on the strata of ${\cal X}^{\pm}\otimes_kk'$ and the closed embeddings 
of ${\cal X}'\otimes_kk'$ into ${\cal X}\otimes_kk'$ definitely do not 
descend to $k$, the structure in $\ell$-adic cohomology that these 
data induce should descend.  We cannot prove this, but we show that 
this approach works at least in the simple case of $U(1,1)$.
\vskip 2mm

In conclusion we wish to thank J.-L.~Waldspurger who communicated to 
us his computation of the sign of the Langlands-Shelstad transfer 
factor which allowed us to make the comparison with the sign factor 
which arises from our geometric approach.  This comparison had been 
requested by R.P.~Langlands at the Princeton conference in his honour 
in October 1996 when a preliminary version of this paper was 
presented.

The results of this paper were obtained during the visits of the first 
author at the Universities of Wuppertal and of K\"oln and the visits 
of the second author at Orsay and at the Institut Emile Borel.  The 
second author especially wishes to thank the members of the 
d\'epartement de math\'ematique de l'Universit\'e de Paris-Sud for 
inviting him and for making his stay a great pleasure.  We both thank 
the Deutsche Forschungsgemeinschaft for its support, as well as the 
European Network (TMR) ``Arithmetic Algebraic Geometry''.
\vskip 5 mm

\centerline{\bf Contents}
\vskip 5 mm

\noindent Part I\par
\vskip 3mm

1. Hermitian forms\par
2. Statement of the Langlands-Shelstad conjecture\par
3. Orbital integrals as numbers of rational points of $k$-schemes\par
4. Statement of the main theorem\par
\vskip 3mm

\noindent Part II\par
\vskip 3mm
5. Stratifications of the $k'$-schemes $X^{\pm}$\par
6. The vector bundle structure\par
7. Proof of the main theorem\par
\vskip 3mm

\noindent Part III\par
\vskip 3mm

8. Descent from $k'$ to $k$\par
9. $U(1,1)$\par
10. Remarks and examples\par
\vskip 3mm

\noindent References\par
\vfill\eject

\centerline{\douzebf PART I}
\vskip 10mm

\centerline{\bf 1. Hermitian forms}
\vskip 5mm

Let $F$ be a non archimedean local field of equal characteristic $p>2$ 
and let $F'$ be an unramified quadratic extension of $F$.  We denote 
by ${\cal O}_F$ and ${\cal O}_{F'}$ the rings of integers of $F$ and 
$F'$ and we fix a uniformizing parameter $\varpi_F$ of $F$ and 
therefore of $F'$.  We denote by $\sigma$ the non trivial element of 
the Galois group ${\rm Gal}(F'/F) \cong {\rm Gal}(k'/k)$ where $k$ and 
$k'$ are the residue fields of $F$ and $F'$.

We have the isomorphism 
$$
{\Bbb Z}/2{\Bbb Z} \buildrel\sim\over\longrightarrow F^\times/{\rm 
N}_{F'/F} F'^\times,\ \varepsilon\mapsto \varpi_F^\varepsilon {\rm 
N}_{F'/F} F'^\times
$$
where ${\rm N}_{F'/F}:F'^\times\rightarrow F^\times$ is the norm map.  
Therefore, for each finite dimensional $F'$-vector space $V$ and each 
$\varepsilon =0,1$ there exists a non degenerate hermitian form on $V$ 
with discriminant $\varpi_F^\varepsilon {\rm N}_{F'/F}F'^\times$.  Any 
two such forms are equivalent.

Let $E$ be a totally ramified separable finite extension of $F$.  The 
tensor product $E'=E\otimes_{F}F'$ is an unramified field extension of 
$E$.  We denote by ${\cal O}_E$ and ${\cal O}_{E'}$ the rings of 
integers of $E$ and $E'$ and we denote again by $\sigma$ the 
automorphism $1\otimes\sigma$ of $E'$.  We denote by $\delta_{E/F}$ 
the exponent of the discriminant ideal ${\frak d}_{{\cal O}_E/{\cal 
O}_F}\subset {\cal O}_F$ of ${\cal O}_E$ over ${\cal O}_F$.

For each $\alpha\in E^\times$ we define a non degenerate hermitian 
form $\Phi_{(\alpha )}$ on the $F'$-vector space $E'$ by
$$
\Phi_{(\alpha )}(e_{1}',e_{2}')={\rm tr}_{E'/F'} (\alpha e_{1}'^\sigma 
e_{2}').
$$
Its discriminant is equal to ${\rm N}_{E/F} (\alpha 
)\varpi_F^{\delta_{E/F}}{\rm N}_{F'/F}F'^\times$ and the dual of the 
${\cal O}_{F'}$-lattice ${\cal O}_{E'}$ with respect to that hermitian
form
$$
({\cal O}_{E'})^{\perp_{(\alpha )}}\buildrel{\rm 
dfn}\over{=\!=}\{e'\in E'\mid \Phi_{(\alpha )} (e',e'')\in {\cal 
O}_{F'},\ \forall e''\in {\cal O}_{E'}\},
$$
is equal to $\alpha^{-1}\varpi_E^{-\delta_{E/F}}{\cal O}_{E'}$ where 
$\varpi_E$ is a uniformizing parameter of $E$ (cf.  [Se] Ch.  III, \S 
3).  As the norm map ${\rm N}_{E/F}: E^\times\rightarrow F^\times$ 
induces an isomorphism from $E^\times /{\cal O}_E^\times$ onto 
$F^\times/ {\cal O}_F^\times$, for any $\alpha^{+}\in 
\varpi_E^{-\delta_{E/F}}{\cal O}_{E'}^\times$ (resp.  $\alpha^{-}\in 
\varpi_E^{1-\delta_{E/F}}{\cal O}_{E'}^\times$) the hermitian form 
$\Phi_{(\alpha^{+})}$ (resp.  $\Phi_{(\alpha^{-})}$) has discriminant $1$ 
(resp.  $\varpi_F$) modulo ${\rm N}_{F'/F}F'^\times$ and the dual of 
the ${\cal O}_{F'}$-lattice ${\cal O}_{E'}$ with respect to that 
hermitian form is equal to ${\cal O}_{E'}$ (resp.  $\varpi_{E}^{-1}{\cal 
O}_{E'}$).  We fix once for all such an element $\alpha^{+}$ (resp.  
$\alpha^{-}$) and we set
$$
\Phi_{E'}^{\pm}\buildrel{\rm dfn}\over{=\!=}\Phi_{(\alpha^{\pm})}.
$$
\vskip 30mm

\centerline{\bf 2. Statement of the Langlands-Shelstad conjecture}
\vskip 5mm

We closely follow [Ko] (in this paper Kottwitz has proved Conjecture 
$2.2$ below in the particular case $n_{1}=1$ and $n_{2}=2$ but for an 
arbitrary local field $F$).

We fix two totally ramified separable finite extensions $E_{1}$ and 
$E_{2}$ of $F$ of degrees $n_{1}$ and $n_{2}$.

We denote by $E_{1}'$ and $E_{2}'$ the unramified quadratic field 
extensions $E_{1}F'$ and $E_{2}F'$ of $E_{1}$ and $E_{2}$.  We denote 
by ${\cal O}_{E_{1}}$, ${\cal O}_{E_{2}}$, ${\cal O}_{E_{1}'}$ and 
${\cal O}_{E_{2}'}$ the rings of integers of $E_{1}$, $E_{2}$, 
$E_{1}'$ and $E_{2}'$.  We fix uniformizing parameters 
$\varpi_{E_{1}}$ and $\varpi_{E_{2}}$ of $E_{1}$ and $E_{2}$ and 
therefore of $E_{1}'$ and $E_{2}'$.

We set $E'=E_{1}'\oplus E_{2}'$.  It is a $F'$-vector space of 
dimension $n_{1}+n_{2}$.  We endow $E'$ with the non degenerate 
hermitian forms
$$
\Phi^{+}=\Phi_{E_{1}'}^{+}\oplus \Phi_{E_{2}'}^{+}
$$
and
$$
\Phi^{-}=\Phi_{E_{1}'}^{-}\oplus \Phi_{E_{2}'}^{-}.
$$
These two forms are equivalent as their discriminants are $1$ and 
$\varpi_F^2$ modulo ${\rm N}_{F'/F}F'^\times$.  Therefore we can find 
$g\in GL_{F'}(E')$ such that
$$
\Phi^{-}(e',e'')=\Phi^{+}(ge',ge'')\qquad (\forall e',e''\in E').
$$

We fix $\gamma_{1}\in E_{1}'^\times$ and $\gamma_{2}\in E_{2}'^\times$ 
such that $\gamma_{1}\gamma_{1}^\sigma= 
\gamma_{2}\gamma_{2}^\sigma=1$, so that $\gamma_{i}$ is a unit in the 
ring ${\cal O}_{E_{i}'}$.  We assume that $E_{i}'= F'[\gamma_{i}]$, 
i.e.  the minimal polynomial $P_{i}(T)\in F'[T]$ of $\gamma_{i}$ has 
degree $n_{i}$.  We assume moreover that the polynomials $P_{1}(T)$ 
and $P_{2}(T)$ are separable and prime with respect to each other.  
Then the diagonal element $(\gamma_{1},\gamma_{2})\in GL_{F'}(E')$ may 
be simultaneously viewed as an elliptic regular semisimple element 
$\gamma^{+}$ in the unitary group
$$
G(F)\buildrel{\rm dfn}\over{=\!=} 
U(E',\Phi^{+})=gU(E',\Phi^{-})g^{-1}\subset GL_{F'}(E'),
$$
as an elliptic regular semisimple element $\gamma^{-}$ in the unitary 
group
$$
U(E',\Phi^{-})\subset GL_{F'}(E')
$$
and as an elliptic $(G,H)$-regular semisimple element $\delta$ in the 
endoscopic group
$$
H(F)=U(E_{1}',\Psi_{1})\times U(E_{2}',\Psi_{2})\subset GL_{F'}(E')
$$
of $G(F)$ where we have set
$$
\Psi_{i}=\Phi_{E_{i}'}^{+}.
$$
The elements $\gamma^{+}$ and $g\gamma^{-}g^{-1}$ of $G(F)$ are 
conjugate in $GL_{F'}(E')$ but are not conjugate in $G(F)$.  The 
conjugacy class of $\delta$ in $H(F)$ is equal to its stable conjugacy 
class (an element of $U(E_{i}',\Psi_{i})\subset GL_{F'}(E_{i}')$ is 
stably conjugate to $\gamma_{i}$ if and only if it has the same 
minimal polynomial as $\gamma_{i}$).

Let $K$ be the hyperspecial maximal compact open subgroup of $G(F)$ 
which fixes the ${\cal O}_{F'}$-lattice ${\cal O}_{E_{1}'}\oplus {\cal 
O}_{E_{2}'}$ of $E'$ (this lattice is selfdual for the hermitian form 
$\Phi^{+}$) and let $K^H$ be the hyperspecial maximal compact open 
subgroup of $H(F)$ which fixes the same lattice.  The $\kappa$-orbital 
integral $O_{\gamma}^\kappa (1_{K})$ is equal to the difference ([Ko]) 
$$\displaylines{
\qquad O_{\gamma}^\kappa (1_{K})=|\{{\cal L}'\subset E'\mid {\cal 
L}'^{\perp^{+}}={\cal L}'\hbox{ and }(\gamma_{1},\gamma_{2}){\cal 
L}'={\cal L}'\}|
\hfill\cr\hfill
 -|\{{\cal L}'\subset E'\mid {\cal L}'^{\perp^{-}}={\cal L}'\hbox{ and 
 }(\gamma_{1},\gamma_{2}){\cal L}'={\cal L}'\}|\qquad}
$$
where the ${\cal L}'$'s are ${\cal O}_{F'}$-lattices and where $(\cdot 
)^{\perp^{\pm}}$ denotes the duality for such lattices with respect to 
the hermitian form $\Phi^{\pm}$.  Similarly the (stable) orbital 
integral $SO_\delta^H (1_{K^H})$ is equal to the product
$$\displaylines{
\qquad SO_\delta^H(1_{K^H})=|\{{\cal M}_{1}'\subset E_{1}'\mid {\cal 
M}_{1}'^{\perp_{1}} ={\cal M}_{1}'\hbox{ and }\gamma_{1}{\cal 
M}_{1}'={\cal M}_{1}'\}|
\hfill\cr\hfill \times |\{{\cal M}_{2}'\subset E_{2}'\mid {\cal 
M}_{2}'^{\perp_{2}} ={\cal M}_{2}'\hbox{ and }\gamma_{2}{\cal 
M}_{2}'={\cal M}_{2}'\}|.\qquad}
$$
where the ${\cal M}_{i}'$'s are ${\cal O}_{F'}$-lattices and where 
$(\cdot )^{\perp_{i}}$ denotes the duality for such lattices with 
respect to the hermitian form $\Psi_{i}$.

As the polynomials $P_{1}(T),P_{2}(T)\in F'[T]$ are prime with respect 
to each other and have coefficients in ${\cal O}_{F'}$ their resultant 
${\rm Res}(P_{1},P_{2})$ is a non zero element in ${\cal O}_{F'}$.  
Let us denote by
$$
r=r(\gamma_{1},\gamma_{2})\geq 0\leqno{\hbox{{\rm (2.1)}}}
$$
its order.  Let us recall that, up to a sign, we have
$$
{\rm Res}(P_{1},P_{2})=\prod_{k_{1}=0}^{n_{1}-1} 
\prod_{k_{2}=0}^{n_{2}-1} (\gamma_{1}^{(k_{1})} -\gamma_{2}^{(k_{2})})
$$
where $\gamma_{i}=\gamma_{i}^{(0)},\ldots ,\gamma_{i}^{(n_{i}-1)}$ are 
the roots of $P_{i}(T)$ in some algebraic closure of $F'$ containing 
$E_{1}'$ and $E_{2}'$.

\th CONJECTURE 2.2 (Langlands-Shelstad)
\enonce
Under the above hypotheses we have
$$
O_\gamma^\kappa (1_{K})=\varepsilon (\gamma_{1},\gamma_{2})q^{r} 
SO_\delta^H(1_{K^H})
$$
where $\varepsilon (\gamma_{1},\gamma_{2})$ is some sign and $q$ is 
the number of elements in the residue field $k$.
\endth

\rem Remark
\endrem
Waldspurger has computed the sign of the Langlands-Shelstad transfer 
factor at any regular semisimple element $\gamma\in G(F)$ which is 
close enough to the identity.  In particular he has proved that
$$
\varepsilon (\gamma_{1},\gamma_{2})=(-1)^{r}
$$
(private communication).  This is exactly the sign which arises from 
our geometric approach (cf.  Theorem 4.2 and \S 8).
\vskip 5mm

\centerline{\bf 3.  Orbital integrals as numbers of rational points of 
$k$-schemes}
\vskip 5mm

Our first goal is to introduce $k$-schemes ${\cal X}^{+}$, ${\cal 
X}^{-}$ and ${\cal Y}_{i}$ ($i=1,2$) such that
$$
O_\gamma^\kappa (1_{K})=|{\cal X}^{+}(k)|-|{\cal X}^{-}(k)|
$$
and
$$
SO_\delta^H(1_{K^H})=|{\cal Y}_{1}(k)|\cdot |{\cal Y}_{2}(k)|.
$$

\th DEFINITION 3.1
\enonce
Let $R$ be a commutative $k'$-algebra.  A $(R\otimes_{k'}{\cal 
O}_{F'})$-lattice in a finite dimensional $F'$-vector space $V$ is a 
$(R\otimes_{k'}{\cal O}_{F'})$-submodule $L$ of $R\otimes_{k'}V$ such 
that there exist ${\cal O}_{F'}$-lattices $L_{0}\subset L_{1}$ in $V$ 
having the following properties:

\decale{\rm (i)} $R\otimes_{k'}L_{0}\subset L\subset 
R\otimes_{k'}L_{1}$,

\decale{\rm (ii)} the $R$-module $L/(R\otimes_{k'}L_{0})$ is locally a 
direct factor of the free $R$-module 
$(R\otimes_{k'}L_{1})/(R\otimes_{k'}L_{0})$.
\vskip 3mm

If $R$ is a commutative $k$-algebra, a $(R\otimes_{k}{\cal 
O}_{F'})$-lattice ${\cal L}$ in a finite dimensional $F'$-vector space 
$V$ is by definition a $((R\otimes_{k}k')\otimes_{k'}{\cal 
O}_{F'})$-lattice in $V$.
\endth

If $R$ is a $k'$-algebra and if $L_{1}$ and $L_{2}$ are 
two $(R\otimes_{k'}{\cal O}_{F'})$-lattices in a finite dimensional 
$F'$-vector space $V$ we set
$$
[L_{1}:L_{2}]={\rm rk}_{R}(L_{1}/(R\otimes_{k'}L_{3})) -{\rm 
rk}_{R}(L_{1}/(R\otimes_{k'}L_{3}))
$$
where $L_{3}$ is any ${\cal O}_{F'}$-lattice in $V$ such that 
$R\otimes_{k'}L_{3}$ is contained in both lattices $L_{1}$ and 
$L_{2}$.  It is a locally constant function on ${\rm Spec}(R)$ with 
integral values.  If $R$ is a $k$-algebra and if ${\cal L}_{1}$ and 
${\cal L}_{2}$ are two $(R\otimes_{k}{\cal O}_{F'})$-lattices we 
define $[{\cal L}_{1}:{\cal L}_{2}]$ by simply replacing $R$ by 
$R\otimes_{k}k'$ in the above definition.

If $V$ is a finite dimensional $F'$-vector space which is equipped 
with a non degenerate hermitian form $\Phi$ and if $R$ is a 
commutative $k$-algebra, we define the dual $(R\otimes_{k}{\cal 
O}_{F'})$-lattice ${\cal L}^{\perp}$ of a $(R\otimes_{k}{\cal 
O}_{F'})$-lattice ${\cal L}$ in the obvious way.  In particular, if we 
have $R\otimes_{k}{\cal L}_{0}\subset {\cal L}\subset 
R\otimes_{k}{\cal L}_{1}$ as in Definition 3.1 we have
$$
R\otimes_{k}{\cal L}_{1}^{\perp}\subset {\cal L}^{\perp}\subset 
R\otimes_{k}{\cal L}_{0}^{\perp}
$$
and, if we identify $(R\otimes_{k}{\cal 
L}_{0}^{\perp})/(R\otimes_{k}{\cal L}_{1}^{\perp})$ with the dual of 
the free $(R\otimes_{k}k')$-module $(R\otimes_{k}{\cal 
L}_{1})/(R\otimes_{k}{\cal L}_{0})$,
$$
{\cal L}^{\perp}/(R\otimes_{k}{\cal L}_{1}^{\perp})\subset 
(R\otimes_{k}{\cal L}_{0}^{\perp})/(R\otimes_{k}{\cal L}_{1}^{\perp})
$$
is the orthogonal of the ${\cal L}/(R\otimes_{k}{\cal L}_{0})\subset 
(R\otimes_{k}{\cal L}_{1})/(R\otimes_{k}{\cal L}_{0})$.
\vskip 3mm

If $R$ is a commutative $k'$-algebra and if $L$ (resp.  $M_{i}$) is a 
$(R\otimes_{k'}{\cal O}_{F'})$-lattice in $E'$ (resp.  $E_{i}'$) we 
define the index of $L$ (resp.  $M_{i}$) as the locally constant 
function
$$
{\rm ind}(L)=[L:R\otimes_{k'}({\cal O}_{E_{1}'}\oplus {\cal 
O}_{E_{2}'})]:{\rm Spec}(R)\rightarrow {\Bbb Z}
$$
(resp.
$$
{\rm ind}(M_{i})=[M_{i}:R\otimes_{k'}{\cal O}_{E_{i}'}]:{\rm 
Spec}(R)\rightarrow {\Bbb Z}\,).
$$
If $R$ is a commutative $k$-algebra and if ${\cal L}$ (resp.  ${\cal 
M}_{i}$) is a $(R\otimes_{k}{\cal O}_{F'})$-lattice in $E'$ (resp.  
$E_{i}'$) we define the index of ${\cal L}$ (resp.  ${\cal M}_{i}$) by 
replacing $R$ by $R\otimes_{k}k'$ in the above definitions.  We then 
have
$$
{\rm ind}({\cal L}^{\perp^{+}})=-{\rm ind}({\cal L})
$$
and
$$
{\rm ind}({\cal L}^{\perp^{-}})=-{\rm ind}({\cal L})+2
$$
where $\perp^{\pm}$ is the duality for lattices with respect to the 
hermitian form $\Phi^{\pm}$ (resp.
$$
{\rm ind}({\cal M}_{i}^{\perp_{i}})=-{\rm ind}({\cal M}_{i})
$$
where $\perp_{i}$ is the duality for lattices with respect to the 
hermitian form $\Psi_{i}$).
\vskip 3mm

For each commutative $k$-algebra $R$ we now set
$$
{\cal X}^{\pm}(R)=\{{\cal L}^{\pm}\mid ({\cal L}^{\pm})^{\perp^{\pm}}= 
{\cal L}^{\pm}\hbox{ and }(1\otimes (\gamma_{1},\gamma_{2})){\cal 
L}^{\pm}={\cal L}^{\pm}\}\leqno{\hbox{{\rm (3.2.1)}}}
$$
and
$$
{\cal Y}_{i}(R)=\{{\cal M}_{i}\mid {\cal M}_{i}^{\perp_{i}}= {\cal 
M}_{i}\hbox{ and } (1\otimes\gamma_{i}){\cal M}_{i}={\cal 
M}_{i}\}\leqno{\hbox{{\rm (3.2.2)}}}
$$
where ${\cal L}^{\pm}$ and ${\cal M}_{i}$ are $(R\otimes_{k}{\cal 
O}_{F'})$-lattices in $E'$ and $E_{i}'$.  If $\varphi :S\rightarrow R$ 
is a homomorphism of commutative $k$-algebras we have obvious base 
change maps ${\cal X}^{\pm}(R)\rightarrow {\cal X}^{\pm}(S)$ and 
${\cal Y}_{i}(R)\rightarrow {\cal Y}_{i}(S)$.

\th PROPOSITION 3.3
\enonce
{\rm (i)} The functor ${\cal X}^{\pm}$ is representable by a $k$-scheme 
which is locally of finite type.

\decale{\rm (ii)} The functor ${\cal Y}_{i}$ is representable by a 
projective $k$-scheme.
\endth

Let us recall that to give a quasi-projective $k$-scheme ${\cal S}$ is 
the same as to give the quasi-projective $k'$-scheme $k'\otimes_k{\cal 
S}$ together with the endomorphism $k'\otimes_{k}{\rm Frob}_{\cal S}$ 
where ${\rm Frob}_{\cal S}$ is the Frobenius endomorphism of ${\cal 
S}$ with respect to $k$ (cf.  [SGA 1](VIII, 7.6)).  Therefore 
Proposition 3.3 is an immediate consequence of the following 
proposition:

\th PROPOSITION 3.4
\enonce
{\rm (i)} The functor $k'\otimes_{k}{\cal X}^{\pm}$ {\rm (}i.e.  the 
restriction of the functor ${\cal X}^{\pm}$ to $k'$-algebras{\rm )} is 
representable by a $k'$-scheme which is an increasing union of 
$(k'\otimes_{k}{\rm Frob}_{{\cal X}^{\pm}})$-stable quasi-projective 
open subsets.

\decale{\rm (ii)} The functor $k'\otimes_{k}{\cal Y}_{i}$ {\rm (}i.e.  
the restriction of the functor ${\cal Y}_{i}$ to $k'$-algebras{\rm )} 
is representable by a projective $k'$-scheme.
\endth

\rem Remark 
\endrem
This proposition is a particular case of a result of D.  Kazhdan and G.  
Lusztig (cf.  [Ka-Lu] \S 2).
\vskip 3mm

Before proving Proposition 3.4 let us give a description of the 
functors $k'\otimes_{k}{\cal X}^{\pm}$ and $k'\otimes_{k}{\cal Y}_{i}$ 
which does not involve the hermitian forms $\Phi^{\pm}$ and $\Psi_{i}$.
\vskip 3mm

For every finite dimensional $F'$-vector space $V$ which is equipped 
with a non degenerate hermitian form $\Phi$ we may split the 
$(k'\otimes_{k}F)$-vector space $k'\otimes_kV$ into
$$
\widetilde V\oplus\widetilde {\widetilde V}
$$
where
$$
\widetilde V=\{x\in k'\otimes_{k}V\mid (1\otimes\alpha')x= 
(\alpha'\otimes 1)x,\ \forall \alpha'\in k'\}
$$
and
$$
\widetilde {\widetilde V}=\{x\in k'\otimes_{k}V\mid (1\otimes\alpha')x= 
(\alpha'^\sigma\otimes 1)x,\ \forall \alpha'\in k'\}.
$$
The $(k'\otimes_{k}F)$-bilinear form $k'\otimes_{k}\Phi$ is then given by
$$
(k'\otimes_{k}\Phi)(\widetilde x_{1}\oplus \widetilde {\widetilde x}_{1}, 
\widetilde x_{2}\oplus \widetilde {\widetilde x}_{2})=\widetilde \Phi 
(\widetilde {\widetilde x}_{1},\widetilde x_{2})+\bigl(\widetilde \Phi 
(\widetilde {\widetilde x}_{2},\widetilde x_{1})\bigr)^{1\otimes\sigma}
$$
for some non degenerate $(k'\otimes_{k}F)$-bilinear form
$$
\widetilde \Phi :\widetilde {\widetilde V}\times \widetilde V 
\rightarrow \{x\in k'\otimes_{k}F'\mid (1\otimes\alpha')x= 
(\alpha'\otimes 1)x,\ \forall \alpha'\in k'\}.
$$
The map $\sigma\otimes_{k}{\rm Id}_{V}: k'\otimes_{k}V\rightarrow 
k'\otimes_{k}V$ induces $(\sigma\otimes_{k}{\rm Id}_{F'})$-linear 
bijections $F:\widetilde V\rightarrow \widetilde {\widetilde V}$ and 
$G:\widetilde {\widetilde V}\rightarrow\widetilde V$.  The maps 
$G\circ F$ and $F\circ G$ are the Frobenius endomorphisms with respect 
to $F'$, i.e.  the identities of the $F'$-vector spaces $\widetilde V$ 
and $\widetilde {\widetilde V}$.  One easily checks that
$$
\widetilde \Phi (F(\widetilde x),G(\widetilde {\widetilde x}))=\bigl( 
\widetilde \Phi (\widetilde {\widetilde x},\widetilde x)\bigr)^{\sigma 
\otimes\sigma}
$$
for every $\widetilde x\in\widetilde V$ and $\widetilde {\widetilde 
x}\in \widetilde {\widetilde V}$.

We may identify $k'\otimes_{k}F$ and $\{x\in k'\otimes_{k}F'\mid 
(1\otimes \alpha')x= (\alpha'\otimes 1)x,\ \forall\alpha'\in k'\}$ 
with $F'$ by
$$
\alpha'\otimes 1\mapsto \alpha'\hbox{ and }\sum_{i} \alpha_{i}'\otimes 
b_{i}'\mapsto\sum_{i}\alpha_{i}'b_{i}'.
$$
Similarly we will identify $\widetilde V$ and $\widetilde {\widetilde 
V}$ with the $F'$-vector spaces $V$ and $k'\otimes_{\sigma ,k'}V$ by
$$
\sum_{i}\alpha_{i}'\otimes v_{i}\mapsto\sum_{i}\alpha_{i}'v_{i}\hbox{ 
and }\sum_{i} \alpha_{i}'\otimes v_{i}\mapsto 
\sum_{i}\alpha_{i}'^\sigma v_{i}.
$$
Then $\Phi : (k'\otimes_{\sigma ,k'}V)\times V\rightarrow F'$ is a non 
degenerate $F'$-bilinear form and we may identify $k'\otimes_{\sigma 
,k'} V$ with the $F'$-linear dual $V^{\ast}$ of $V$.  We have 
$\sigma$-linear bijections $F:V\rightarrow V^{\ast}$ and $G: 
V^{\ast}\rightarrow V$.  The maps $G\circ F$ and $F\circ G$ are the 
identities of the $F'$-vector spaces $V$ and $V^{\ast}$.  We have
$$
\langle F(x),G(x^{\ast} )\rangle =\sigma\bigl(\langle x^{\ast} 
,x\rangle\bigr)
$$
for every $x\in V$ and $x^{\ast}\in V^{\ast}$.

\th LEMMA 3.5
\enonce
Let us denote by $\perp$ the duality {\rm (}for lattices{\rm )} with 
respect to $\Phi$ and let us assume moreover that $V$ admits a 
selfdual ${\cal O}_{F'}$-lattice ${\cal L}_{0}$.

Then, for every commutative $k'$-algebra $R$ there is a natural 
bijection between the set of $(R\otimes_{k}{\cal O}_{F'})$-lattices 
${\cal L}$ in $V$ such that ${\cal L}^{\perp}={\cal L}$ and the set of 
$(R\otimes_{k'}{\cal O}_{F'})$-lattices $L$ in $V$ such that 
$[L:R\otimes_{k'}{\cal L}_{0}]=0$.
\endth

\rem Proof
\endrem
The relation between ${\cal L}$ and $L$ is
$$
{\cal L}=L\oplus L^\perp\subset (R\otimes_{k'}V)\oplus (R\otimes_{k'} 
V^{\ast} ) =R\otimes_{k}V
$$
where $L^\perp$ is the dual $(R\otimes_{k'}F')$-lattice of $L$ in 
$R\otimes_{k'}V^{\ast}$.
\hfill\hfill$\square$
\vskip 3mm

If $u$ is a unitary automorphism of $(V,\Phi)$ the 
$(k'\otimes_kF')$-linear automorphism $1\otimes u$ of $k'\otimes_kV$ 
is the direct sum $u\oplus ({}^{\rm t}u)^{-1}$ where ${}^{\rm t}u$ is 
the transposed endomorphism of $u$.  Moreover we have the relations
$$
{}^{\rm t}u\circ F\circ u=F,\quad  u\circ G\circ {}^{\rm t}u=G. 
$$

In particular, we can take $(V,\Phi )=(E_{i}',\Phi_{E_{i}'}^{\pm})$, 
$(V,\Phi )=(E',\Phi^{\pm})$ or $(V,\Phi )= (E_{i}',\Psi_{i})$ and we 
get the $\sigma$-linear bijections $F_{i}^{\pm} :E_{i}'\rightarrow 
E_{i}'^{\ast}$, $G_{i}^{\pm}:E_{i}'^{\ast}\rightarrow E_{i}'$, 
$F^{\pm} =F_{1}^{\pm}\oplus F_{2}^{\pm}:E'\rightarrow E'^{\ast}$, 
$G^{\pm} =G_{1}^{\pm}\oplus G_{2}^{\pm}:E'^{\ast}\rightarrow E'$, 
$F_{i}=F_{i}^{+}:E_{i}'\rightarrow E_{i}'^{\ast}$ and 
$G_{i}=G_{i}^{+}: E_{i}'^{\ast}\rightarrow E_{i}'$.
\vskip 3mm

Now, for any commutative $k'$-algebra $R$ we set
$$
X^{\pm}(R)=\{L^{\pm}\mid {\rm ind}(L^{\pm})\equiv\delta^{\pm} \hbox{ 
and }(1\otimes (\gamma_{1},\gamma_{2}))L^{\pm}= 
L^{\pm}\}\leqno{\hbox{{\rm (3.6.1)}}}
$$
where $\delta^{+}=0$ and $\delta^{-}=1$ and where the $L^{\pm}$'s are 
$(R\otimes_{k'}{\cal O}_{F'})$-lattices in $E'$, and we set
$$
Y_{i}(R)=\{M_{i}\mid {\rm ind}(M_{i})\equiv 0\hbox{ and } 
(1\otimes\gamma_{i})M_{i}=M_{i}\}\leqno{\hbox{{\rm (3.6.2)}}}
$$
where the $M_{i}$'s are $(R\otimes_{k'}{\cal O}_{F'})$-lattices in 
$E_{i}'$.  If $\varphi :S\rightarrow R$ is a homomorphism of 
commutative $k$-algebras we have obvious base change maps 
$X^{\pm}(R)\rightarrow X^{\pm}(S)$ and $Y_{i}(R)\rightarrow Y_{i}(S)$.  
We denote by
$$
F_{X^{\pm}}:X^{\pm}\rightarrow X^{\pm}\leqno{\hbox{{\rm (3.6.3)}}}
$$
and
$$
F_{Y_{i}}:Y_{i}\rightarrow Y_{i}\leqno{\hbox{{\rm (3.6.4)}}}
$$
the functor endomorphisms which are given on the $R$-valued points by
$$
L^{\pm}\mapsto \bigl(F_R^{\pm}(L^{\pm})\bigr)^{\perp^{\pm}}= 
G_R^{\pm}\bigl((L^{\pm})^{\perp^{\pm}}\bigr)
$$
and
$$
M_{i}\mapsto \bigl(F_{i,R}(M_{i})\bigr)^{\perp_{i}}= 
G_{i,R}\bigl((M_{i})^{\perp_{i}}\bigr)\,)
$$
where $F_R^{\pm} : R\otimes_{k'}E'\rightarrow (R\otimes_{k'}E' 
)^{\ast}$, $G_R^{\pm} : (R\otimes_{k'}E')^{\ast}\rightarrow 
R\otimes_{k'}E'$, $F_{i,R} : R\otimes_{k'}E_{i}'\rightarrow 
(R\otimes_{k'}E_{i}')^{\ast}$ and $G_{i,R} : 
(R\otimes_{k'}E_{i}')^{\ast}\rightarrow R\otimes_{k'}E_{i}'$) are the 
natural $\sigma$-linear extensions of $F^{\pm}$, $G^{\pm}$, $F_{i}$ 
and $G_{i}$.

Then it follows from the above discussion that

\th LEMMA 3.7
\enonce
We have natural identifications
$$
(k'\otimes_{k}{\cal X}^{\pm},k'\otimes_{k}{\rm Frob}_{{\cal 
X}^{\pm}})=(X^{\pm},F_{X^{\pm}})
$$
and
$$
(k'\otimes_{k}{\cal Y}_{i},k'\otimes_{k}{\rm Frob}_{{\cal 
Y}_{i}})=(Y_{i},F_{Y_{i}}).
$$
\hfill\hfill$\square$
\endth

Let
$$
m_{i}=m(\gamma_{i})\geq 0\leqno{\hbox{{\rm (3.8)}}}
$$
be the conductor of ${\cal O}_{F'}[\gamma_{i}]$ in ${\cal 
O}_{E_{i}'}$, i.e.  the smallest non negative integer $m$ such that
$$
\varpi_{E_{i}}^m{\cal O}_{E_{i}'}\subset {\cal O}_{F'}[\gamma_{i}]\subset 
{\cal O}_{E_{i}'}.
$$

The proof of Proposition 3.4 is based on the next two lemmas:

\th LEMMA 3.9
\enonce
Let $R$ be a commutative $k'$-algebra and let $L_{i}$ be a 
$(R\otimes_{k'}{\cal O}_{F'})$-lattice in $E_{i}'$ such that 
$\gamma_{i}L_{i}\subset L_{i}$ with constant index.  Then there exists 
an integer $\ell$ such that
$$
{\rm ind}(L_{i})\leq \ell\leq {\rm ind}(L_{i})+m_{i}
$$
and
$$
\varpi_{E_{i}}^{m_{i}-\ell}(R\otimes_{k'}{\cal O}_{E_{i}'})\subset 
L_{i}\subset \varpi_{E_{i}}^{-\ell}(R\otimes_{k'}{\cal O}_{E_{i}'}).
$$

In particular we automatically have
$$
\varpi_{E_{i}}^{m_{i}-{\rm ind}(L_{i})}(R\otimes_{k'}{\cal 
O}_{E_{i}'})\subset L_{i}\subset \varpi_{E_{i}}^{-m_{i}-{\rm 
ind}(L_{i})}(R\otimes_{k'}{\cal O}_{E_{i}'}).
$$
\endth

\rem Proof 
\endrem
We have
$$
\varpi_{E_{i}}^{m_{i}}{\cal O}_{E_{i}'}L_{i}\subset {\cal 
O}_{F'}[\gamma_{i}]L_{i}\subset L_{i}\subset {\cal O}_{E_{i}'}L_{i}
$$
and we only need to check that ${\cal O}_{E_{i}'}L_{i}$ is equal to 
$\varpi_{E_{i}}^{-\ell}(R\otimes_{k'}{\cal O}_{E_{i}'})$ for some 
integer $\ell$.

But, as $L_{i}$ is a $(R\otimes_{k'}{\cal O}_{F'})$-lattice there 
exist ${\cal O}_{F'}$-lattices ${\cal L}_{0}\subset {\cal L}_{1}$ in 
$E_{i}'$ such that $R\otimes_{k'}{\cal L}_{0}\subset L_{i}\subset 
R\otimes_{k'}{\cal L}_{1}$ and we have
$$
R\otimes_{k'}({\cal O}_{E_{i}'}{\cal L}_{0})\subset {\cal 
O}_{E_{i}'}L_{i}\subset R\otimes_{k'}({\cal O}_{E_{i}'}{\cal L}_{1})
$$
with ${\cal O}_{E_{i}'}{\cal L}_{0}=\varpi_{E_{i}}^{-\ell_{0}}{\cal 
O}_{E_{i}'}$ and ${\cal O}_{E_{i}'}{\cal L}_{1}= 
\varpi_{E_{i}}^{-\ell_{1}}{\cal O}_{E_{i}'}$ for some integers 
$\ell_{0}\leq \ell_{1}$. Now the multiplication by $\varpi_{E_{i}}$ 
on the free $R$-module
$$
\varpi_{E_{i}}^{-\ell_{1}}(R\otimes_{k'}{\cal O}_{E_{i}'})/ 
\varpi_{E_{i}}^{-\ell_{0}}(R\otimes_{k'}{\cal O}_{E_{i}'})
$$
is a regular nilpotent endomorphism and the $R$-submodule
$$
{\cal 
O}_{E_{i}'}L_{i}/\varpi_{E_{i}}^{-\ell_{0}}(R\otimes_{k'}{\cal 
O}_{E_{i}'})\subset \varpi_{E_{i}}^{-\ell_{1}}(R\otimes_{k'}{\cal 
O}_{E_{i}'})/\varpi_{E_{i}}^{-\ell_{0}}(R\otimes_{k'}{\cal O}_{E_{i}'})
$$
is locally a direct factor.  Therefore this $R$-submodule is equal to 
the $R$-submodule
$$\varpi_{E_{i}}^{-\ell}(R\otimes_{k'}{\cal O}_{E_{i}'})/ 
\varpi_{E_{i}}^{-\ell_{0}}(R\otimes_{k'}{\cal O}_{E_{i}'})
$$
for some integer $\ell$ with $\ell_{0}\leq \ell\leq \ell_{1}$ and the 
proof of the lemma is complete.
\hfill\hfill$\square$
\vskip 3mm

For each point $L^{\pm}$ in $X^{\pm}$ with value in some field extension 
$K$ of $k'$ let us denote by $B_{i}^{\pm}$ the intersection of $L^{\pm}$ 
with
$K\otimes_{k'}E_{i}'\subset K\otimes_{k'}E'$ and by $C_{i}^{\pm}$ the 
projection of $L^{\pm}$ on $K\otimes_{k'}E_{i}'$.  Then
$$
B_{i}^{\pm}\subset C_{i}^{\pm}\subset 
K\otimes_{k'}E_{i}'\leqno{\hbox{{\rm (3.10.1)}}}
$$
are $(K\otimes_{k'}{\cal O}_{F'})$-lattices in $E_{i}'$ and $L^{\pm}$ is the 
graph of an isomorphism of $(K\otimes_{k'}{\cal O}_{F'})$-modules of 
finite dimension over $K$
$$
\iota^{\pm}:C_{1}^{\pm}/B_{1}^{\pm}\isomorphism 
C_{2}^{\pm}/B_{2}^{\pm}.\leqno{\hbox{{\rm (3.10.2)}}}
$$
Moreover, if we set
$$
b_{i}^{\pm}=b_{i}^{\pm}(L^{+}):={\rm ind}(B_{i}^{\pm})\leq 
c_{i}^{\pm}=c_{i}^{\pm}(L^{+}):={\rm ind}(C_{i}^{\pm})\leqno{\hbox{{\rm (3.10.3)}}}
$$
we have
$$
b_{1}^{\pm}+c_{2}^{\pm}=b_{2}^{\pm}+c_{1}^{\pm}={\rm 
ind}(L^{\pm})=\delta^{\pm}\leqno{\hbox{{\rm (3.10.4)}}}
$$
with $\delta^{+}=0$ and $\delta^{-}=1$.

If $L^{\pm}$ is now a point in $X^{\pm}$ with value in some 
commutative $k'$-algebra $R$ we define as above integers
$$
b_{i}^{\pm}({\frak p})=b_{i}^{\pm}({\rm Frac}(R/{\frak 
p})\otimes_{R}L^{\pm})\leq c_{i}^{\pm}({\frak p})=c_{i}^{\pm}({\rm 
Frac}(R/{\frak p})\otimes_{R}L^{\pm})
$$
for each prime ideal ${\frak p}$ in $R$.  We thus have functions 
$b_{i}^{\pm},c_{i}^{\pm}:{\rm Spec}(R)\rightarrow {\Bbb Z}$ which are 
easily seen to be semicontinuous: for each integer $\lambda_{i}$ the 
set of points in ${\rm Spec}(R)$ such that 
$b_{i}^{\pm}\leq\lambda_{i}$ (or equivalently such that 
$c_{j}^{\pm}\geq\delta^{\pm}-\lambda_{i}$ if $\{i,j\}=\{1,2\}$) is 
open for the Zariski topology.

\th LEMMA 3.11
\enonce
Let $L^{\pm}$ be a point in $X^{\pm}$ with value in some commutative 
$k'$-algebra $R$.  Then the functions $b_{i}^{\pm},c_{i}^{\pm}:{\rm 
Spec}(R)\rightarrow {\Bbb Z}$ satisfy the inequalities
$$
b_{i}^{\pm}\leq c_{i}^{\pm}\leq b_{i}^{\pm}+r.
$$
\endth

\rem Proof 
\endrem
We may assume that $R$ is a field extension of $k'$ and therefore we 
may introduce the lattices $B_{i}^{\pm},C_{i}^{\pm}$ and the 
isomorphism $\iota^{\pm}$.  As $\iota^{\pm}$ exchanges the 
multiplications by $\gamma_{1}$ and $\gamma_{2}$ and as 
$P_{1}(\gamma_{1})=0$ and $P_{2}(\gamma_{2})=0$ we have
$$
P_{2}(\gamma_{1})C_{1}^{\pm}\subset B_{1}^{\pm}
$$
and
$$
P_{1}(\gamma_{2})C_{2}^{\pm}\subset B_{2}^{\pm}.
$$
But $P_{2}(\gamma_{1})$ and $P_{1}(\gamma_{1})$ are of order 
$r=r(\gamma_{1},\gamma_{2})$ in $E_{1}'$ and $E_{2}'$, so that
$$
b_{i}^{\pm}\leq c_{i}^{\pm}\leq b_{i}^{\pm}+r
$$
as required.
\hfill\hfill$\square$
\vskip 3mm

\rem Proof of Proposition $3.4$
\endrem
Let us begin with Part (ii).  It follows from Lemma 3.9 that, for any 
commutative $k'$-algebra $R$, $Y_{i}(R)$ may be identified with the 
set of $(R\otimes_{k'}{\cal O}_{F'})$-lattices $M_{i}$ in $E_{i}'$ such 
that
$$
\varpi_{E_{i}}^{m_{i}}(R\otimes_{k'}{\cal O}_{F'})\subset 
M_{i}\subset\varpi_{E_{i}}^{-m_{i}}(R\otimes_{k'}{\cal O}_{F'}),
$$
$$
{\rm rk}_{R}(M_{i}/\varpi_{E_{i}}^{m_{i}}(R\otimes_{k'}{\cal 
O}_{F'}))=m_{i}
$$
and
$$
(1\otimes\gamma_{i})M_{i}\subset M_{i}.
$$
Therefore, $Y_{i}$ is representable by a closed $k'$-subscheme of the 
Grassmann variety of $m_{i}$-planes in the $2m_{i}$-dimensional 
$k'$-vector space
$$
\varpi_{E_{i}}^{-m_{i}}{\cal O}_{F'}/\varpi_{E_{i}}^{m_{i}}{\cal 
O}_{F'}
$$
and Part (ii) is proved.

The proof of Part (i) is similar to the proof of Part (ii).  For each 
pair of integers $(\lambda_{1},\lambda_{2})$ let us 
consider the open subfunctor
$$
X^{\pm}(\lambda_{1},\lambda_{2})=\{L^{\pm}\in X^{\pm}\mid 
b_{1}^{\pm}\leq\lambda_{1}\hbox{ and 
}b_{2}^{\pm}\leq\lambda_{2}\}.\leqno{\hbox{{\rm (3.12.1)}}}
$$
We only need to prove that this open subfunctor is representable by a 
quasi-projective $k'$-scheme.

It follows from Lemma 3.11 that $X^{\pm}(\lambda_{1},\lambda_{2})$ is 
contained in the closed subfunctor
$$
X^{\pm}[\delta^{\pm}-\lambda_{2}-r,\delta^{\pm}-\lambda_{1}-r]= 
\{L^{\pm}\in X^{\pm}\mid b_{1}^{\pm}\geq \delta^{\pm}-\lambda_{2}-r 
\hbox{ and }b_{2}^{\pm}\geq \delta^{\pm}-\lambda_{1}-r\}
$$
of $X^{\pm}$.  Therefore it is sufficient to prove that, for any pair 
of integers $(\mu_{1},\mu_{2})$ the closed subfunctor
$$
X^{\pm}[\mu_{1},\mu_{2}]=\{L^{\pm}\in X^{\pm}\mid b_{1}^{\pm}\geq 
\mu_{1}\hbox{ and }b_{2}^{\pm}\geq \mu_{2}\}\leqno{\hbox{{\rm (3.12.2)}}}
$$
is representable by a projective $k'$-scheme.  But it follows from 
Lemma 3.9 that the functor $X^{\pm}[\mu_{1},\mu_{2}]$ is a closed 
subfunctor of the functor which associates to any commutative 
$k'$-algebra $R$ the set of $(R\otimes_{k'}{\cal O}_{F'})$-lattices 
$L^{\pm}$ in $E'$ such that
$$\displaylines{
\qquad\varpi_{E_{1}}^{m_{1}-\mu_{1}}(R\otimes_{k'}{\cal 
O}_{E_{1}'})\oplus \varpi_{E_{2}}^{m_{2}-\mu_{2}} (R\otimes_{k'}{\cal 
O}_{E_{2}'})
\hfill\cr\hfill
\subset 
L^{\pm}\subset\varpi_{E_{1}}^{\mu_{2}-m_{1}-\delta^{\pm}}(R\otimes_{k'}{\cal 
O}_{E_{1}'})\oplus \varpi_{E_{2}}^{\mu_{1}-m_{2}-\delta^{\pm}} 
(R\otimes_{k'}{\cal O}_{E_{2}'}),\qquad}
$$
$$
{\rm rk}_{R}\bigl(L^{\pm}/(\varpi_{E_{1}}^{m_{1}-\mu_{1} 
}(R\otimes_{k'}{\cal O}_{E_{1}'})\oplus \varpi_{E_{2}}^{m_{2}-\mu_{2}} 
(R\otimes_{k'}{\cal O}_{E_{2}'}))\bigr)= 
m_{1}+m_{2}-\mu_{1}-\mu_{2}+\delta^{\pm}
$$
and
$$
(1\otimes (\gamma_{1},\gamma_{2}))L^{\pm}\subset L^{\pm}.
$$
Therefore, if we set $N=m_{1}+m_{2}-\mu_{1}-\mu_{2}+\delta^{\pm}$ the 
functor $X^{\pm}[\mu_{1},\mu_{2}]$ is representable by a closed 
$k'$-subscheme of the Grassmann variety of $N$-planes in the 
$2N$-dimensional $k'$-vector space
$$
(\varpi_{E_{1}}^{m_{1}-\mu_{1}}{\cal O}_{E_{1}'}\oplus 
\varpi_{E_{2}}^{m_{2}-\mu_{2}} {\cal O}_{E_{2}'})/ 
(\varpi_{E_{1}}^{\mu_{2}-m_{1}-\delta^{\pm}}{\cal O}_{E_{1}'}\oplus 
\varpi_{E_{2}}^{\mu_{1}-m_{2}-\delta^{\pm}} {\cal O}_{E_{2}'})
$$
and Part (ii) is proved.
\hfill\hfill$\square$
\vskip 3mm

\rem Remark $3.13$
\endrem
In fact the two $k'$-schemes $X^{+}$ and $X^{-}$ are isomorphic.  More 
precisely, there are two $k'$-scheme isomorphisms
$$
\alpha_{1},\alpha_{2}:X^{-}\isomorphism X^{+}
$$
given by
$$
\alpha_{1}(L^{-})=(\varpi_{E_{1}}\oplus 1)L^{-}
$$
and
$$
\alpha_{2}(L^{-})=(1\oplus \varpi_{E_{2}})L^{-}.
$$
Moreover the squares of $k'$-scheme morphisms
$$\diagram{
X^{-}&\kern -1mm\longmaprightover{F_{X^{-}}}{10mm}\kern -1mm&X^{-}\cr
\mapdownleft{\alpha_{1}}&&\mapdownright{\alpha_{2}}\cr
X^{+}&\kern -1mm\longmaprightunder{F_{X^{+}}}{10mm}\kern -1mm&X^{+}\cr}
\quad\hbox{ and }\quad\diagram{
X^{-}&\kern -1mm\longmaprightover{F_{X^{-}}}{10mm}\kern -1mm&X^{-}\cr
\mapdownleft{\alpha_{2}}&&\mapdownright{\alpha_{1}}\cr
X^{+}&\kern -1mm\longmaprightunder{F_{X^{+}}}{10mm}\kern -1mm&X^{+}\cr}
$$
commute.
\vskip 3mm

\rem Remark $3.14$
\endrem
Let $G_{i}$ be the Grassmann variety of $m_{i}$-planes in the 
$2m_{i}$-dimensional $k'$-vector space
$$
W_{i}=\varpi_{E_{i}}^{-m_{i}}{\cal O}_{E_{i}'}/\varpi_{E_{i}}^{m_{i}}{\cal 
O}_{E_{i}'}
$$
and let $N_{i}$ be the regular nilpotent endomorphism of $W_{i}$ which 
is induced by the multiplication by $\varpi_{E_{i}}$.  In the course 
of the proof of Proposition 3.4 we have identified $Y_{i}$ with a 
closed subscheme of $G_{i}$.

Let $Z_{i}$ be the closed subscheme of $G_{i}$ of $m_{i}$-planes 
$D_{i}\subset W_{i}$ such that
$$
N_{i}^{m_{i}+j}(D_{i})\subset D_{i}
$$
for every non negative integer $j$. Then, as $\varpi_{E_{i}}^{m_{i}}{\cal 
O}_{E_{i}'}\subset {\cal O}_{F'}[\gamma_{i}]$ by definition of the 
conductor $m_{i}$, we have
$$
Y_{i}\subset Z_{i}.
$$

The structure of $Z_{i}$ is much simpler than that of $Y_{i}$.  In 
fact, it follows from Lemma 3.15 below that either $m_{i}=0$ and 
$Z_{i}=G_{i}={\rm Spec}(k')$, or $m_{i}>0$ and
$$
Z_{i}=\bigcup_{j=1}^{m_{i}-1}Z_{i,j}
$$
where
$$
Z_{i,j}=\{D_{i}\subset W_{i}\mid {\rm dim}(D_{i})=m_{i}\hbox{ and }{\rm 
Ker}(N_{i}^{j})\subset D_{i}\subset {\rm Ker}(N_{i}^{m_{i}+j})\}.
$$
Moreover, if $m_{i}>0$ the irreducible components of $Z_{i}$ are 
exactly the $Z_{i,j}$ for $j=1,\ldots ,m_{i}-1$ and the intersection 
of any two of these irreducible components is again a Grassmann 
variety.
\vskip 3mm

\th LEMMA 3.15
\enonce 
Let $E$ be a finite dimensional $k'$-vector space which is equipped 
with a regular nilpotent endomorphism $N$.  We denote by $e$ the 
dimension of $E$.  Let $m$ be a non negative integer and let $F$ be a 
$k'$-vector subspace of $E$ which is $N^{m+j}$-stable for every $j\geq 
0$.  Let $\ell$ be a positive integer.

If $F$ is not contained in ${\rm Im}\bigl(N^\ell\bigr)$ it 
automatically contains ${\rm Im}\bigl(N^{m+\ell -1}\bigr)$ and 
therefore its dimension is at least $e-m-\ell +1$ {\rm (}and even at 
least $e-m-\ell +2$ if $m\geq 1$ as ${\rm Im}\bigl(N^{m+\ell 
-1}\bigr)\subset{\rm Im}\bigl(N^\ell\bigr)$ in this case{\rm )}.

If $F$ does not contain ${\rm Ker}\bigl(N^\ell\bigr)$ it is 
automatically contained in ${\rm Ker}\bigl(N^{m +\ell -1}\bigr)$ and 
therefore its dimension is at most $m+\ell -1$ {\rm (}and even at most 
$m+\ell -2$ if $m\geq 1$ as ${\rm Ker}\bigl(N^{m +\ell 
-1}\bigr)\supset {\rm Ker}\bigl(N^\ell\bigr)$ in this case{\rm )}.
\endth

\rem Proof
\endrem
If $F\not\subset {\rm Im}\bigl(N^\ell\bigr)={\rm Ker}\bigl(N^{e-\ell 
}\bigr)$ there exists $x\in F$ such that $N^{e-\ell}(x)\not= 0$.  Let 
$j$ be the unique non negative integer such that $N^{e-\ell 
+j}(x)\not= 0$ and $N^{e-\ell +j+1}(x)= 0$.  Then we may consider the 
vector subspace of $E$ generated by 
$N^{m+j}(x),N^{m+j+1}(x),\ldots,\break N^{e-\ell +j}(x)$.  It is 
contained in $F$ by hypothesis and it is easily seen to be equal to 
${\rm Ker}\bigl(N^{e-m-\ell +1}\bigr)$.  The first assertion of the 
lemma follows.

The second assertion of the lemma is another formulation of the first 
one.
\hfill\hfill$\square$
\vskip 5mm

\centerline{\bf 4.  Statement of the main theorem}
\vskip 5mm

On the $k'$-scheme $X^{\pm}$ we have a translation automorphism 
$\tau^{\pm}:X^{\pm}\rightarrow X^{\pm}$ which is given by
$$
\tau^{\pm}(L^{\pm})=(\varpi_{E_{1}}\oplus\varpi_{E_{2}}^{-1})L^{\pm}
$$
and which satisfies the anticommutation relation
$$
\tau^{\pm}\circ F_{X^{\pm}}=F_{X^{\pm}}\circ (\tau^{\pm})^{-1}.
$$
In fact, with the notation of Remark 3.13 we have
$$
\tau^{+}=\alpha_{1}\circ \alpha_{2}^{-1}
$$
and
$$
\tau^{-}=\alpha_{2}^{-1}\circ \alpha_{1}.
$$
In particular $\tau^{-}$ and $\tau^{+}$ are exchanged by the 
isomorphisms $\alpha_{1}$ and $\alpha_{2}$.

\rem Remark 
\endrem
The action of ${\Bbb Z}$ on $X^{\pm}$ generated by the translation 
automorphism $\tau^{\pm}$ has been introduced by D.  Kazhdan and G.  
Lusztig in [Ka-Lu] \S 2.
\vskip 3mm

{\it Let us first assume that $r=2r'$ is even}.  Then the projective closed 
subscheme (cf.  (3.12.2))
$$\eqalign{
X^{+}[-r',-r']=&\{L^{+}\in X^{+}\mid b_{1}^{+}\geq -r'\hbox{ and 
}b_{2}^{+}\geq -r'\}\cr
=&\{L^{+}\in X^{+}\mid -r'\leq b_{i}^{+}\leq c_{i}^{+}\leq r'\}\cr}
$$
may be viewed as a ``fundamental domain'' for the translation 
$\tau^{+}$ on the $k'$-scheme $X^{+}$: for each integer $n$ we have
$$\eqalign{
(\tau^{+})^{n}X^{+}[-r',-r']=&X^{+}[-r'-n,-r'+n]\cr
=&\{L^{+}\in X^{+}\mid -r'-n\leq b_{1}^{+}\leq c_{1}^{+}\leq 
r'-n\}\cr
=&\{L^{+}\in X^{+}\mid -r'+n\leq b_{2}^{+}\leq c_{2}^{+}\leq 
r'+n\},\cr}
$$
so that
$$
X^{+}=\bigcup_{n\in{\Bbb Z}}(\tau^{+})^{n}X^{+}[-r',-r']
$$
and the intersection
$$
X^{+}[-r',-r']\cap (\tau^{+})^{n}X^{+}[-r',-r']=X^{+}[-{\rm 
inf}(r',r'+n),-{\rm inf}(r',r'-n)]
$$
is empty if $|n|>r$. Moreover, as
$$
b_{i}^{+}(F_{X^{+}}L^{+})=-c_{i}^{+}(L^{+})
$$
all the fixed points by $F_{X^{+}}$ in $X^{+}$ are contained in 
$X^{+}[-r',-r']$ and we have
$$
|{\cal X}^{+}(k)|=\big|(X^{+}[-r',-r'])^{F_{X^{+}}}\big|.  
\leqno{\hbox{{\rm (4.1.1)}}}
$$
Similarly, as
$$
b_{i}^{-}(F_{X^{-}}L^{-})=1-c_{i}^{-}(L^{-})
$$
all the fixed points by $F_{X^{-}}$ in $X^{-}$ are contained in 
$$
X^{-}[1-r',1-r']=\{L^{-}\in X^{-}\mid 1-r'\leq b_{i}^{-}\leq 
c_{i}^{-}\leq r'\}
$$
and we have
$$
|{\cal X}^{-}(k)|=\big|(X^{-}[1-r',1-r'])^{F_{X^{-}}}\big|.  
\leqno{\hbox{{\rm (4.1.2)}}}
$$
\vskip 1mm

{\it Let us now assume that $r=2r'+1$ is odd}.  Then the projective closed 
subscheme (cf.  (3.12.2))
$$\eqalign{
X^{-}[-r',-r']=&\{L^{-}\in X^{-}\mid b_{1}^{-}\geq -r'\hbox{ and 
}b_{2}^{-}\geq -r'\}\cr
=&\{L^{-}\in X^{-}\mid -r'\leq b_{i}^{-}\leq c_{i}^{-}\leq 1+r'\}\cr}
$$
may be viewed as a ``fundamental domain'' for the translation 
$\tau^{-}$ on the $k'$-scheme $X^{-}$: for each integer $n$ we have
$$\eqalign{
(\tau^{-})^{n}X^{-}[-r',-r']=&X^{-}[-r'-n,-r'+n]\cr
=&\{L^{-}\in X^{-}\mid -r'-n\leq b_{1}^{-}\leq c_{1}^{-}\leq 
1+r'-n\}\cr
=&\{L^{-}\in X^{-}\mid -r'+n\leq b_{2}^{-}\leq c_{2}^{-}\leq 
1+r'+n\},\cr}
$$
so that
$$
X^{-}=\bigcup_{n\in{\Bbb Z}}(\tau^{-})^{n}X^{-}[-r',-r']
$$
and the intersection
$$
X^{-}[-r',-r']\cap (\tau^{-})^{n}X^{-}[-r',-r']=X^{-}[-{\rm 
inf}(r',r'+n),-{\rm inf}(r',r'-n)]
$$
is empty if $|n|>r$. Moreover, as
$$
b_{i}^{-}(F_{X^{-}}L^{-})=1-c_{i}^{-}(L^{-})
$$
all the fixed points by $F_{X^{-}}$ in $X^{-}$ are contained in 
$X^{-}[-r',-r']$ and we have
$$
|{\cal X}^{-}(k)|=\big|(X^{-}[-r',-r'])^{F_{X^{-}}}\big|.  
\leqno{\hbox{{\rm (4.1.3)}}}
$$
Similarly, as
$$
b_{i}^{+}(F_{X^{+}}L^{+})=-c_{i}^{+}(L^{+})
$$
all the fixed points by $F_{X^{+}}$ in $X^{+}$ are contained in 
$$
X^{+}[-r',-r']=\{L^{+}\in X^{+}\mid -r'\leq b_{i}^{+}\leq c_{i}^{+}\leq 
r'\}
$$
and we have
$$
|{\cal X}^{+}(k)|=\big|(X^{+}[-r',-r'])^{F_{X^{+}}}\big|.  
\leqno{\hbox{{\rm (4.1.4)}}}
$$

\th THEOREM 4.2
\enonce
Let us set
$$\eqalign{
(X,F_{X})&=\bigl(X^{+}[-r',-r'],F_{X^{+}}\big|X^{+}[-r',-r']\bigr)\cr
(X',F_{X'})&=\bigl(X^{-}[1-r',1-r'],F_{X^{-}}\big|X^{-}[1-r',1-r']\bigr)\cr}
$$
if $r=2r'$ is even and
$$\eqalign{
(X,F_{X})&=\bigl(X^{-}[-r',-r'],F_{X^{-}}\big|X^{-}[-r',-r']\bigr)\cr
(X',F_{X'})&=\bigl(X^{+}[-r',-r'],F_{X^{+}}\big|X^{+}[-r',-r']\bigr)\cr}
$$
if $r=2r'+1$ is odd, so that $X$ and $X'$ are projective $k'$-schemes 
and that $F_{X}$ and $F_{X'}$ are the Frobenius endomorphisms for some 
$k$-rational structures on $X$ and $X'$ respectively.

Then, for every even positive integer $f$ we have
$$
\big|X^{F_{X}^{f}}\big|- \big|X'^{F_{X'}^{f}}\big|= q^{fr}\cdot 
\big|Y_{1}^{F_{Y_{1}}^{f}}\big|\cdot\big|Y_{2}^{F_{Y_{1}}^{f}}\big|.
$$
\endth

The proof of Theorem 4.2 will be given in \S 7.
\vskip 3mm

\rem Remark
\endrem
Taking into account the relations 
$$
X^{F_{X}}={\cal X}^{\pm}(k)\quad\hbox{ and }\quad X'^{F_{X'}}={\cal X}^{\mp}(k)
$$
where we have put $\pm =+$ if $r$ is even and $\pm=-$ if $r$ is odd 
(cf.  $(4.1.1)$ to $(4.1.4)$), the Langlands-Shelstad conjecture is 
equivalent to the analogous statement of Theorem $4.2$ for $f=1$.  An 
obvious extrapolation of the Langlands-Shelstad conjecture is that the 
equalities in the statement of Theorem $4.2$ should hold for 
every positive odd integer $f$.
\vfill\eject

\centerline{\douzebf PART II}
\vskip 10mm

\centerline{\bf 5.  Stratifications of the $k'$-scheme $X^{\pm}$}
\vskip 5mm

With the notation of (3.10.1) to (3.10.4) a point $L^{\pm}$ in 
$X^{\pm}$ with values in some field extension $K$ of $k'$ may be given 
as a triple
$$
(B_{1}^{\pm}\subset C_{1}^{\pm}, B_{2}^{\pm}\subset 
C_{2}^{\pm},\iota^{\pm})
$$
where $B_{i}^{\pm}$ and $C_{i}^{\pm}$ are $(K\otimes_{k'}{\cal 
O}_{F'})$-lattices in $E_{i}'$ such that
$$
b_{1}^{\pm}+c_{2}^{\pm}=b_{2}^{\pm}+c_{1}^{\pm}=\delta^{\pm}
$$
with $\delta^{+}=0$ and $\delta^{-}=1$, and where
$$
\iota^{\pm}:C_{1}^{\pm}/B_{1}^{\pm}\isomorphism C_{2}^{\pm}/B_{2}^{\pm}
$$
is an isomorphism of $(K\otimes_{k'}{\cal O}_{F'})$-modules (of finite 
dimension over $K$) which exchanges the automorphisms induced by the 
multiplications by $\gamma_{1}$ and $\gamma_{2}$.  It may also be 
given as a triple
$$
(B_{1}^{\pm}, C_{2}^{\pm},f_{1}^{\pm})
$$
(resp.
$$
(C_{1}^{\pm}, B_{2}^{\pm},f_{2}^{\pm})\,)
$$
where $B_{1}^{\pm}$ and $C_{2}^{\pm}$ (resp. $C_{1}^{\pm}$ and 
$B_{2}^{\pm}$) are as before and where
$$
f_{1}^{\pm}:C_{2}^{\pm}\rightarrow (K\otimes_{k'}E_{1}')/B_{1}^{\pm}
$$
(resp.
$$
f_{2}^{\pm}:C_{1}^{\pm}\rightarrow (K\otimes_{k'}E_{2}')/B_{2}^{\pm}\,)
$$
is a morphism of $(K\otimes_{k'}{\cal O}_{F'})$-modules which exchanges 
the automorphisms induced by the multiplications by $\gamma_{1}$ and 
$\gamma_{2}$.
\vskip 3mm

For $i=1,2$ and each integer $j$ let us denote by
$$
U_{i,j}^{\pm}
$$
the locally closed subset of the  $k'$-scheme $X^{\pm}$ which 
is defined by the condition
$$
b_{i}^{\pm}(L^{\pm})=j.
$$
For $i=1,2$ the family $(U_{i,j}^{\pm})_{j\in {\Bbb Z}}$ of locally 
closed subsets is a partition of $X^{\pm}$.  For each integer $j$ we 
have $k'$-scheme morphisms
$$
\pi_{1,j}^{\pm}:U_{1,j}^{\pm}\rightarrow Y_{1}\times_{k'}Y_{2}
$$
and
$$
\pi_{2,j}^{\pm}:U_{2,j}^{\pm}\rightarrow Y_{1}\times_{k'}Y_{2}
$$
which are defined by
$$
\pi_{1,j}^{\pm}(L^{\pm})=(\varpi_{E_{1}}^{j}B_{1}^{\pm}, 
\varpi_{E_{2}}^{\delta^{\pm}-j}C_{2}^{\pm})
$$
and
$$
\pi_{2,j}^{\pm}(L^{\pm})=(\varpi_{E_{1}}^{\delta^{\pm}-j}C_{1}^{\pm}, 
\varpi_{E_{2}}^{j}B_{2}^{\pm}).
$$

\rem Remark 
\endrem
The above $k'$-scheme morphisms are nothing else than the restriction 
of the map $q_{N}$ introduced by Kazhdan and Lusztig in [Ka-Lu] \S 5.
\vskip 3mm

It is obvious that
$$
(\tau^{\pm})^{k}(U_{i,j}^{\pm})=U_{i,j-k}^{\pm}
$$
where $\tau^{\pm}$ is the translation automorphism introduced in 
Section 4 and that
$$
\pi_{i,j}^{\pm}=\pi_{i,j-k}^{\pm}\circ (\tau^{\pm})^{k}.
$$
It is also obvious that
$$
U_{i,j}^{+}=\alpha_{i}(U_{i,j+1}^{-})
$$
where $\alpha_{i}:X^{-}\isomorphism X^{+}$ is the isomorphism 
introduced in Remark 3.13 and that
$$
\pi_{i,j}^{+}\circ\alpha_{i}=\pi_{i,j+1}^{-}.
$$
Moreover we have
$$
F_{X^{\pm}}(U_{1,j}^{\pm})=U_{2,j}^{\pm}\quad\hbox{ and }\quad 
F_{X^{\pm}}(U_{2,j}^{\pm})=U_{1,j}^{\pm}
$$
and if we simply denote by $F_{U}$ the $k'$-morphisms 
$U_{1,j}^{\pm}\rightarrow U_{2,j}^{\pm}$ and $U_{2,j}^{\pm}\rightarrow 
U_{1,j}^{\pm}$ which are induced by $F_{X^{\pm}}$, the squares of 
$k'$-scheme morphisms
$$\diagram{
U_{1,j}^{\pm}&\kern -5.5mm\longmaprightover{F_{U}}{19mm}\kern 
-5.5mm&U_{2,j}^{\pm}\cr
\mapdownleft{\pi_{1,j}^{+}}&&\mapdownright{\pi_{2,j}^{+}}\cr
Y_{1}\times_{k'}Y_{2}&\kern -1mm\longmaprightunder{F_{Y_{1}}\times 
F_{Y_{2}}}{10mm}\kern -1mm&Y_{1}\times_{k'}Y_{2}\cr}
\quad\hbox{ and }\quad\diagram{
U_{2,j}^{\pm}&\kern -5.5mm\longmaprightover{F_{U}}{19mm}\kern 
-5.5mm&U_{1,j}^{\pm}\cr
\mapdownleft{\pi_{2,j}^{+}}&&\mapdownright{\pi_{1,j}^{+}}\cr
Y_{1}\times_{k'}Y_{2}&\kern -1mm\longmaprightunder{F_{Y_{1}}\times 
F_{Y_{2}}}{10mm}\kern -1mm&Y_{1}\times_{k'}Y_{2}\cr}
$$
commute.
\vskip 20mm

\centerline{\bf 6. The vector bundle structure}
\vskip 5mm

\th THEOREM 6.1
\enonce
For $i=1,2$ and each $j\in {\Bbb Z}$ the morphism
$$
\pi_{i,j}^{\pm}: U_{i,j}^{\pm}\rightarrow Y_{1}\times_{k'}Y_{2}
$$
is a rank $r$ vector bundle.
\endth

\rem Proof 
\endrem
First of all it is not difficult to see that $\pi_{i,j}^{\pm}$ has a 
natural structure of generalized vector bundle structure in the sense 
of Grothendieck, the fiber of which through a $K$-rational point 
$L^{\pm}\in X^{\pm}$ is the vector space
$$
{\rm Hom}_{{\cal O}_{F'},\gamma}(C_{2}^{\pm}, 
(K\otimes_{k'}E_{1}')/B_{1}^{\pm})
$$
of ${\cal O}_{F'}$-linear maps $f_{1}^{\pm}:C_{2}^{\pm} \rightarrow 
(K\otimes_{k'}E_{1}')/B_{1}^{\pm}$ which exchange the automorphisms 
induced by the multiplications by $\gamma_{2}$ and $\gamma_{1}$ if 
$i=1$, and the vector space
$$
{\rm Hom}_{{\cal 
O}_{F'},\gamma}(C_{1}^{\pm},(K\otimes_{k'}E_{2}')/B_{2}^{\pm})
$$
of ${\cal O}_{F'}$-linear maps $f_{2}^{\pm}:C_{1}^{\pm} \rightarrow 
(K\otimes_{k'}E_{2}')/B_{2}^{\pm}$ which exchange the automorphisms 
induced by the multiplications by $\gamma_{1}$ and $\gamma_{2}$ if 
$i=2$.

Therefore it suffices to show that the rank of this generalized vector 
bundle is constant of value $r$. But this is a direct consequence of 
the following proposition.
\hfill\hfill$\square$
\vskip 3mm

\th PROPOSITION 6.2
\enonce
For $i=1,2$ let $M_{i}$ be a $(K\otimes_{k'}{\cal O}_{F'})$-lattice in 
$E_{i}'$ where $K$ is some field extension of $k'$.

We consider ${\cal O}_{F'}[\gamma_{1}]$ and ${\cal 
O}_{F'}[\gamma_{2}]$ as quotients of the $2$-dimensional regular local 
ring ${\cal O}_{F'}[[T]]$ by sending $T$ onto $\gamma_{1}$ and 
$\gamma_{2}$ respectively.
 
Then the dimension of the $K$-vector space
$$
{\rm Hom}_{K\otimes_{k'}{\cal 
O}_{F'}[[T]]}(M_{2},(K\otimes_{k'}E_{1}')/M_{1})
$$
is equal to $r$.
\endth

\rem Proof
\endrem
For simplicity we will only consider the case $K=k'$, the proof for 
$K$ arbitrary being the same.  In order to prove the proposition it is 
sufficient to prove that the complex
$$
R{\rm Hom}_{{\cal O}_{F'}[[T]]}(M_{2},E_{1}'/M_{1})
$$
is concentrated in degree $0$ and that its Euler-Poincar\'e 
characteristic is equal to $r$.

\decale{\rm (a)} We have a distinguished triangle
$$
R{\rm Hom}_{{\cal O}_{F'}[[T]]}(M_{2},M_{1})\rightarrow R{\rm 
Hom}_{{\cal O}_{F'}[[T]]}(M_{2},E_{1}')\rightarrow R{\rm Hom}_{{\cal 
O}_{F'}[[T]]}(M_{2},E_{1}'/M_{1})\rightarrow
$$
in the derived category of $k'$-vector spaces.  Moreover $R{\rm 
Hom}_{{\cal O}_{F'}[[T]]}(M_{2},E_{1}')$ is zero in this 
derived category as the endomorphism
$$
R{\rm Hom}_{{\cal O}_{F'}[[T]]}(P_{2}(\gamma_{2}),E_{1}')=R{\rm 
Hom}_{{\cal O}_{F'}[[T]]}(M_{2},P_{2}(\gamma_{1}))
$$
of this complex is at the same time zero and an isomorphism.  We thus 
have an isomorphism
$$
R{\rm Hom}_{{\cal O}_{F'}[[T]]}(M_{2},E_{1}'/M_{1})\isomorphism R{\rm 
Hom}_{{\cal O}_{F'}[[T]]}(M_{2},M_{1})[1]
$$
in the derived category.

The complex
$$
R{\rm Hom}_{{\cal O}_{F'}[[T]]}(M_{2},M_{1})
$$
is concentrated in degrees $1$ and $2$ as ${\rm Hom}_{{\cal 
O}_{F'}[[T]]}(M_{2},M_{1})$ is obviously zero and as ${\cal 
O}_{F'}[[T]]$ is a regular local ring of dimension $2$.  To prove the 
proposition it is thus sufficient to prove that
$$
{\rm Ext}_{{\cal O}_{F'}[[T]]}^2(M_{2},M_{1})=(0)
$$
and that the Euler-Poincar\'e characteristic of $R{\rm Hom}_{{\cal 
O}_{F'}[[T]]}(M_{2},M_{1})$ is equal to $-r$.

\decale{\rm (b)} Let us take a more geometric point of view.  Let $S$ 
be the germ of surface ${\rm Spec}({\cal O}_{F'}[[T]])$ and, for 
$\alpha =1,2$ let
$$
\iota_\alpha :C_\alpha ={\rm Spec}({\cal O}_{F'}[\gamma_\alpha 
])\hookrightarrow S
$$
be the germ of curve with equation $P_\alpha (T)=0$.  By hypothesis 
these two curves are integral and the support of their intersection is 
the origin $s=(\varpi_{F}=0,T=0)$ of $S$.  For $\alpha =1,2$ the 
${\cal O}_{F'}[\gamma_\alpha ]$-module $M_\alpha$ defines a torsion 
free coherent ${\cal O}_{C_\alpha}$-module of generic rank $1$ that we 
will also denote by $M_\alpha$.  Clearly we have
$$
R{\rm Hom}_{{\cal O}_{F'}[[T]]}(M_{2},M_{1})=R{\rm Hom}_{{\cal 
O}_S}(\iota_{2,\ast}M_{2},\iota_{1,\ast}M_{1})
$$
and the support of $R{\cal H}{\it om}_{{\cal 
O}_S}(\iota_{2,\ast}M_{2},\iota_{1,\ast}M_{1})$ is $\{s\}$.

\decale{\rm (c)} Let us now show that
$$
{\rm Ext}_{{\cal O}_S}^2(\iota_{2,\ast}M_{2},\iota_{1,\ast}M_{1})=(0).
$$
Let $\widetilde{M}_{2}$ be the saturated module
$$
\widetilde{M}_{2}={\cal O}_{E_{2}'}M_{2}\subset {\cal O}_{E_{2}'}
$$
From the geometric point of view $\widetilde{M}_{2}$ is nothing 
else than the free ${\cal O}_{\widetilde{C}_{2}}$-module of rank one 
$\pi_{2}^{\ast} M_{2}$ where $\pi_{2}:\widetilde{C}_{2}\rightarrow 
C_{2}$ is the normalization of the curve $C_{2}$.  We have an 
injective map of ${\cal O}_S$-modules
$$
\iota_{2,\ast}M_{2}\hookrightarrow \iota_{2,\ast}\pi_{2,\ast} 
\widetilde{M}_{2}
$$
and then a surjective map of $k'$-vector spaces
$$
{\rm Ext}_{{\cal 
O}_S}^2(\iota_{2,\ast}\pi_{2,\ast}\widetilde{M}_{2},\iota_{1,\ast}M_{1}) 
\twoheadrightarrow {\rm Ext}_{{\cal 
O}_S}^2(\iota_{2,\ast}M_{2},\iota_{1,\ast}M_{1})
$$
(all the ${\rm Ext}_{{\cal O}_S}^3$'s are zero as $S$ is a regular 
surface).  It is thus sufficient to prove that
$$
{\rm Ext}_{{\cal O}_S}^2(\iota_{2,\ast}\pi_{2,\ast}\widetilde{M}_{2}, 
\iota_{1,\ast}M_{1})=(0).
$$
But, by Grothendieck duality (cf.  [Ha]) this last Ext group is 
isomorphic to
$$
{\rm Ext}_{{\cal O}_{\widetilde C_{2}}}^2(\widetilde{M}_{2}, 
L(\iota_{2}\circ\pi_{2})^!\iota_{1,\ast}M_{1})
$$
and as $\widetilde{M}_{2}$ is a free $\widetilde{{\cal 
O}}_{\widetilde{C}_{2}}$-module of rank $1$ it is isomorphic to
$$
{\cal H}^2L(\iota_{2}\circ\pi_{2})^!\iota_{1,\ast}M_{1}.
$$
As $S$ is a regular surface and $C_{2}$ is a Cartier divisor on $S$, 
the dualizing complexes $K_S$ and $K_{\widetilde{C}_{1}}$ are of the 
form $\omega_S[2]$ and $\omega_{\widetilde{C}_{2}}[1]$ where 
$\omega_S$ and $\omega_{\widetilde{C}_{2}}$ are invertible modules and 
$L(\iota_{2}\circ\pi_{2})^!\iota_{1,\ast}M_{1}$ is isomorphic to 
$L(\iota_{2}\circ\pi_{2})^{\ast}\iota_{1,\ast}M_{1}[-1]$.  As this 
last complex is obviously concentrated in degrees $\leq 1$ the 
assertion follows.

\decale{\rm (d)} Finally, thanks to a result of Deligne (see 
Proposition $6.3$ and its corollary below) it follows that the 
Euler-Poincar\'e characteristic of $R{\rm Hom}_{{\cal 
O}_S}(\iota_{2,\ast}M_{2},\iota_{1,\ast}M_{1})$ is equal
to
$$
(-1)^{{\rm dim}(C_{1})}{\rm mult}(C_{2},C_{1})
$$
where ${\rm mult}(C_{2},C_{1})$ is the intersection multiplicity of 
$C_{2}$ and $C_{1}$ on $S$.  By definition this intersection 
multiplicity is equal to $r$.  This concludes the proof of Proposition 
$6.2$.
\hfill\hfill$\square$
\vskip 3mm

\th PROPOSITION  6.3 (Deligne)
\enonce
Let $(X,x)$ be a germ of smooth variety over a field $k$.  For 
simplicity we assume that $x$ is rational over $k$.  Let $\iota_{1}: 
Y_{1}\hookrightarrow X$ and $\iota_{2}: Y_{2}\hookrightarrow X$ be two 
integral closed subschemes of $X$ such that ${\rm dim}(Y_{1})+{\rm 
dim}(Y_{2})={\rm dim}(X)$ and that $(Y_{1}\cap Y_{2})_{\rm 
red}=\{x\}$.  Let $K_{1}$ and $K_{2}$ be two bounded complexes of flat 
quasi-coherent ${\cal O}_X$-modules.  We assume that all the 
cohomology sheaves of $K_{1}$ {\rm (}resp.  $K_{2}${\rm )} are 
coherent ${\cal O}_X$-modules with support in $Y_{1}$ {\rm (}resp.  
$Y_{2}${\rm )}.

Then $\Gamma (X,K_{1}\otimes_{{\cal O}_X}K_{2})$ is a bounded complex 
of $k$-vector spaces with finite dimensional cohomology and its 
Euler-Poincar\'e characteristic
$$
\sum_n(-1)^n{\rm dim}_{k}\Gamma (X,{\cal H}^n(K_{1}\otimes_{{\cal 
O}_X}K_{2}))
$$
is equal to
$$
m(Y_{1},Y_{2}){\rm rk}_{\eta_{1}}(K_{1}){\rm rk}_{\eta_{2}}(K_{2})
$$
where $m(Y_{1},Y_{2})$ is the intersection multiplicity of $Y_{1}$ and 
$Y_{2}$, i.e.
$$
m(Y_{1},Y_{2})=\sum_n(-1)^n{\rm dim}_{k}{\rm Tor}_n^{{\cal 
O}_{X,x}}({\cal O}_{Y_{1},x},{\cal O}_{Y_{2},x}),
$$
and where, for $\alpha =1,2$
$$
{\rm rk}_{\eta_\alpha }(K_\alpha )=\sum_n(-1)^n{\rm dim}_{\kappa 
(\eta_\alpha )}{\cal H}^n(K_\alpha )_{\eta_\alpha}
$$
is the rank of $K_\alpha$ at the generic point $\eta_\alpha$ of 
$Y_\alpha$.
\endth

\rem Proof
\endrem
See [SGA4${1\over 2}$] [Cycle] Thm.  2.3.8 (iii).  Deligne's argument 
can be easily modified to avoid the use of $\ell$-adic cohomology.
\hfill\hfill$\square$
\vskip 3mm

\th COROLLARY 6.4
\enonce
Let us assume moreover that $Y_{2}$ is a complete intersection in $X$.  
Let ${\cal M}_{1}$ and ${\cal M}_{2}$ be two coherent modules on 
$Y_{1}$ and $Y_{2}$ respectively with generic rank $r_{1}$ and 
$r_{2}$.  Then the complex of $k$-vector spaces
$$
R{\rm Hom}_{{\cal O}_X}(\iota_{2,\ast}{\cal M}_{2},\iota_{1,\ast}{\cal 
M}_{1})
$$
has bounded and finite dimensional cohomology and its Euler-Poincar\'e 
characteristic is equal to
$$
(-1)^{{\rm dim}(Y_{1})}m(Y_{2},Y_{1})r_{1}r_{2}.
$$
\endth

\rem Proof
\endrem
For $\alpha =1,2$ let
$$
K_\alpha\rightarrow\iota_{\alpha ,\ast}{\cal M}_\alpha\rightarrow 0
$$
be a finite resolution of the ${\cal O}_X$-module $\iota_{\alpha 
,\ast}{\cal M}_\alpha$ by free ${\cal O}_X$-modules of finite rank.  
Obviously the complexes $K_{1}$ and $K_{2}$ satisfy the hypotheses of 
the above proposition and ${\rm rk}_{\eta_\alpha}(K_\alpha 
)=r_\alpha$.

In the derived category of $k$-vector spaces the complex
$$
R{\rm Hom}_{{\cal O}_X}(\iota_{2,\ast}{\cal M}_{2},\iota_{1,\ast}{\cal 
M}_{1})
$$
is isomorphic to
$$
{\rm Hom}_{{\cal O}_X}(K_{2},K_{1})=\Gamma (X,K_{2}^\vee\otimes_{{\cal 
O}_X}K_{1})
$$
where $K_{2}^\vee ={\rm Hom}_{{\cal O}_X}(K_{2},{\cal O}_X)$ is also a 
bounded complex of free ${\cal O}_X$-modules of finite rank.  As $X$ 
is local and smooth over $k$, ${\cal O}_X[{\rm dim}(X)]$ is a 
dualizing complex for $X$ and as $Y_{2}$ is a complete intersection in 
$X$, ${\cal O}_{Y_{2}}[{\rm dim}(Y_{2})]$ is a dualizing complex for 
$Y_{2}$.  Therefore, by Grothendieck duality we have
$$
R{\cal H}{\it om}_{{\cal O}_X}(\iota_{2,\ast}{\cal M}_{2},{\cal 
O}_X[{\rm dim}(X)])= \iota_{2,\ast}R{\cal H}{\it om}_{{\cal 
O}_{Y_{2}}}({\cal M}_{2},{\cal O}_{Y_{2}}[{\rm dim}(Y_{2})])
$$
and thus $K_{2}^\vee$ is isomorphic to $\iota_{2,\ast}R{\cal H}{\it 
om}_{{\cal O}_{Y_{2}}}({\cal M}_{2},{\cal O}_{Y_{2}})[-{\rm 
dim}(Y_{1})]$ in the derived category of ${\cal O}_X$-modules.  In 
particular
$$
{\rm rk}_{\eta_{2}}(K_{2})=(-1)^{-{\rm dim}(Y_{1})}r_{2}.
$$

Now the corollary is an immediate consequence of the proposition.
\hfill\hfill\cqfd
\vskip 5mm

\centerline{\bf 7. Proof of the main theorem}
\vskip 5mm

Let us define two closed embeddings
$$
i_{1},i_{2}:X'\,\longhookrightarrowover{}{6mm}\,X\leqno{\hbox{{\rm 
(7.1)}}}
$$
by
$$
i_{1}(L^{-})=\alpha_{2}(L^{-})=(1\oplus \varpi_{E_{2}})L^{-}
$$
and
$$
i_{2}(L^{-})=\alpha_{1}(L^{-})=(\varpi_{E_{1}}\oplus 1)L^{-}
$$
if $r$ is even and by
$$
i_{1}(L^{+})=\alpha_{1}^{-1}(L^{+})=(\varpi_{E_{1}}^{-1}\oplus 1)L^{+}
$$
and
$$
i_{2}(L^{+})=\alpha_{2}^{-1}(L^{+})=(1\oplus \varpi_{E_{2}}^{-1})L^{+}
$$
if $r$ is odd.

By definition we have
$$
U_{1}:=X-i_{1}(X')=U_{1,-r'}^{+}\subset X\subset X^{+}
$$
and
$$
U_{2}:=X-i_{2}(X')=U_{2,-r'}^{+}\subset X\subset X^{+}
$$
if $r$ is even and
$$
U_{1}:=X-i_{1}(X')=U_{1,-r'}^{-}\subset X\subset X^{-}
$$
and
$$
U_{2}:=X-i_{2}(X')=U_{2,-r'}^{-}\subset X\subset X^{-}
$$
if $r$ is odd. Therefore we have projections
$$
\pi_{1}:U_{1}\rightarrow Y_{1}\times_{k'}Y_{2}\hbox{ and 
}\pi_{2}:U_{2}\rightarrow Y_{1}\times_{k'}Y_{2}\leqno{\hbox{{\rm 
(7.2)}}}
$$
which are both rank $r$ vector bundles by Theorem 6.1 
($\pi_{i}=\pi_{i,-r'}^{+}$ if $r$ is even and 
$\pi_{i}=\pi_{i,-r'}^{-}$ if $r$ is odd).

Theorem 4.2 follows.
\hfill\hfill$\square$

\vfill\eject

\centerline{\douzebf Part III}
\vskip 10mm

\centerline{\bf 8. Descent from $k'$ to $k$}
\vskip 5mm

Up to this point our approach to the fundamental lemma was essentially 
elementary in that it was based on counting points of algebraic 
varieties.  However to treat the extensions of odd degree of $k$ and 
thus the original statement of Langlands and Shelstad we envisage the 
use of $\ell$-adic cohomology and the Grothendieck fixed point formula 
(cf.  [Gr]).  Indeed whereas all the structures considered so far are 
only defined over $k'$ we predict that their $\ell$-adic cohomology 
descends to $k$.  Even though we cannot prove this assertion we will 
now explain more precisely its meaning.
\vskip 3mm

In \S 3 we have introduce the projective $k'$-schemes $Y_{1}$ and 
$Y_{2}$ together with the Frobenius endomorphisms $F_{Y_{1}}$ and 
$F_{Y_{2}}$ relative to $k$.  We set
$$
(Y,F_{Y})=(Y_{1},F_{Y_{1}})\times_{k'}(Y_{2},F_{Y_{2}}).
$$
In the statement of Theorem 4.2 we have introduced the projective 
$k'$-schemes $X$ and $X'$ together with the Frobenius endomorphisms 
$F_{X}$ and $F_{X'}$ relative to $k$.  In the course of the proof of 
Theorem 4.2 we have introduced the closed embeddings $i_{1}:X' 
\,\longhookrightarrowover{}{6mm}\, X$ and $i_{2}:X' 
\,\longhookrightarrowover{}{6mm}\, X$ and the projections 
$\pi_{1}:U_{1}:=X-i_{1}(X') \,\longhookrightarrowover{}{6mm}\, Y$ and 
$\pi_{2}:U_{2}:=X-i_{2}(X') \,\longhookrightarrowover{}{6mm}\, Y$.  We 
have shown that $\pi_{1}$ and $\pi_{2}$ have a natural structure of 
rank $r$ vector bundle.  Let us denote by
$$
j_{1}:U_{1} \,\longhookrightarrowover{}{6mm}\, X\quad\hbox{ and }\quad 
j_{2}:U_{2} \,\longhookrightarrowover{}{6mm}\, X
$$
the obvious open embeddings.

Then the diagrams
$$\diagram{
X'&\!\longhookrightarrowover{i_{1}}{10mm}\!&X& 
\!\longhookleftarrowover{j_{1}}{10mm}\!&U_{1}\cr
&&&&\mapdownright{\pi_{1}}\cr
&&&&Y\cr
\noalign{\vskip 6mm}\cr
X'&\!\longhookrightarrowover{i_{2}}{10mm}\!&X& 
\!\longhookleftarrowover{j_{2}}{10mm}\!&U_{2}\cr
&&&&\mapdownright{\pi_{2}}\cr
&&&&Y\cr}\leqno{\hbox{{\rm (8.1)}}}
$$
are exchanged by the Frobenius endomorphisms $F_{X'}$, $F_{X}$ and 
$F_{Y}$: there are unique $k'$-scheme 
morphisms $F_{U}:U_{1}\longrightarrow U_{2}$ and 
$F_{U}:U_{2}\longrightarrow U_{1}$ such that the diagrams
$$\diagram{
X'&\!\longhookrightarrowover{i_{1}}{10mm}\!&X& 
\!\longhookleftarrowover{j_{1}}{10mm}\!&U_{1}& 
\!\longtwoheadrightarrowover{\pi_{1}}{10mm}\!&Y\cr
\mapdownleft{F_{X'}}&&\mapdownright{F_{X}}&& 
\mapdownright{F_{U}}&&\mapdownright{F_{Y}}\cr
X'&\!\longhookrightarrowover{i_{2}}{10mm}\!&X& 
\!\longhookleftarrowover{j_{2}}{10mm}\!&U_{2}&
\!\longtwoheadrightarrowover{\pi_{2}}{10mm}\!&Y\cr
\noalign{\vskip 6mm}
X'&\!\longhookrightarrowover{i_{2}}{10mm}\!&X& 
\!\longhookleftarrowover{j_{2}}{10mm}\!&U_{2}&
\!\longtwoheadrightarrowover{\pi_{2}}{10mm}\!&Y\cr
\mapdownleft{F_{X'}}&&\mapdownright{F_{X}}&& 
\mapdownright{F_{U}}&&\mapdownright{F_{Y}}\cr
X'&\!\longhookrightarrowover{i_{1}}{10mm}\!&X& 
\!\longhookleftarrowover{j_{1}}{10mm}\!&U_{1}& 
\!\longtwoheadrightarrowover{\pi_{1}}{10mm}\!&Y\cr}\leqno{\hbox{{\rm 
(8.2)}}}
$$
are commutative.

Let us fix a prime number $\ell$ which is not equal to the 
characteristic of $k$.  For any quasi-projective (resp.  projective) 
variety $Z$ over $k'$ we set
$$
R\Gamma_{\rm c}(Z)=R\Gamma_{\rm c}(\overline k\otimes_{k'}Z,{\Bbb 
Q}_{\ell})
$$
(resp.
$$
R\Gamma (Z)=R\Gamma (\overline k\otimes_{k'}Z,{\Bbb Q}_{\ell})\,).
$$

We consider the following morphisms of distinguished triangles in the 
derived category $D_{\rm c}^{\rm b}({\rm Spec}(k'),{\Bbb Q}_{\ell})$
$$\diagram{
R\Gamma_{\rm c}(U_{2})&\!\longmaprightover{j_{2,!}}{10mm}\!&R\Gamma 
(X)&\!\longmaprightover{i_{2}^{\ast}}{10mm}\!&R\Gamma (X')& 
\!\longmaprightover{\partial_{2}}{10mm}\!&R\Gamma_{\rm c}(U_{2})[1]\cr
\mapdownleft{F_{U}^{\ast}}&&\mapdownleft{F_{X}^{\ast}} 
&&\mapdownright{F_{X'}^{\ast}}&&\mapdownright{F_{U}^{\ast} [1]}\cr
R\Gamma_{\rm c}(U_{1})&\!\longmaprightover{j_{1,!}}{10mm}\!&R\Gamma 
(X)&\!\longmaprightover{i_{1}^{\ast}}{10mm}\!&R\Gamma (X')& 
\!\longmaprightover{\partial_{1}}{10mm}\!&R\Gamma_{\rm 
c}(U_{1})[1]\cr}
$$
and
$$\diagram{
R\Gamma_{\rm c}(U_{1})&\!\longmaprightover{j_{1,!}}{10mm}\!&R\Gamma 
(X)&\!\longmaprightover{i_{1}^{\ast}}{10mm}\!&R\Gamma (X')& 
\!\longmaprightover{\partial_{1}}{10mm}\!&R\Gamma_{\rm c}(U_{1})[1]\cr
\mapdownleft{F_{U}^{\ast}}&&\mapdownleft{F_{X}^{\ast}} 
&&\mapdownright{F_{X'}^{\ast}}&&\mapdownright{F_{U}^{\ast} [1]}\cr
R\Gamma_{\rm c}(U_{2})&\!\longmaprightover{j_{2,!}}{10mm}\!&R\Gamma 
(X)&\!\longmaprightover{i_{2}^{\ast}}{10mm}\!&R\Gamma (X')& 
\!\longmaprightover{\partial_{2}}{10mm}\!&R\Gamma_{\rm 
c}(U_{2})[1]\cr}
$$
(cf.  [SGA4] Exp. XVII, 5.1.16).

We also consider the morphisms 
$$
\pi_{1}^{!}:R\Gamma (Y)[-2r](-r)\longrightarrow R\Gamma_{\rm c}(U_{1})
$$
and
$$
\pi_{2}^{!}:R\Gamma (Y)[-2r](-r)\longrightarrow R\Gamma_{\rm c}(U_{2})
$$
in that derived category which are induced by the trace morphisms
$$
{\rm tr}_{1}:R\pi_{1,!}{\Bbb Q}_{\ell}\longrightarrow {\Bbb 
Q}_{\ell}[-2r](-r)\quad\hbox{ and }\quad{\rm tr}_{2}:R\pi_{2,!}{\Bbb 
Q}_{\ell}\longrightarrow {\Bbb Q}_{\ell}[-2r](-r)
$$
(cf.  [SGA 4] Exp. XVIII, Thm. 2.9).

\th PROPOSITION 8.3
\enonce
The morphisms $\pi_{1}^{!}$ and $\pi_{2}^{!}$ are both isomorphisms 
and the squares in $D_{{\rm c}}^{{\rm b}}({\rm Spec}(k'),{\Bbb 
Q}_{\ell})$ which are induced by the right squares of the diagrams 
$(8.2)$,
$$\diagram{
R\Gamma (Y)[-2r](-r)& \!\longmaprightover{\pi_{2}^{!}}{10mm}\!  
&R\Gamma_{\rm c}(U_{2})\cr
\mapdownleft{F_{Y}^{\ast}[-2r](-r)}&&\mapdownright{F_{U}^{\ast}}\cr
R\Gamma (Y)[-2r](-r)&\!\longmaprightunder{\pi_{1}^{!}}{10mm}\!& 
R\Gamma_{\rm c}(U_{1})\cr}
$$
and
$$\diagram{
R\Gamma (Y)[-2r](-r)&\!\longmaprightover{\pi_{1}^{!}}{10mm}\!& 
R\Gamma_{\rm c}(U_{1})\cr
\mapdownleft{F_{Y}^{\ast}[-2r](-r)}&&\mapdownright{F_{U}^{\ast}}\cr
R\Gamma (Y)[-2r](-r)&\!\longmaprightunder{\pi_{2}^{!}}{10mm}\!& 
R\Gamma_{\rm c}(U_{2})\cr}
$$
are commutative.
\endth

\rem Proof
\endrem
As $\pi_{1}$ and $\pi_{2}$ are both rank $r$ vector bundles (cf.  
Theorem 6.1) the trace morphisms ${\rm tr}_{1}$ and ${\rm tr}_{2}$, 
and therefore the induced morphisms $\pi_{1}^{!}$ and $\pi_{2}^{!}$, 
are all isomorphisms (cf. [SGA 4] Exp. XVIII, Thm. 2.9).

The commutative right squares of the diagrams $(8.2)$ 
induce isomorphisms
$$
f_{1}: F_{Y}^{\ast} R\pi_{2,!}{\Bbb Q}_{\ell}\isomorphism 
R\pi_{1,!}{\Bbb Q}_{\ell}
$$
and
$$
f_{2}: F_{Y}^{\ast} R\pi_{1,!}{\Bbb Q}_{\ell}\isomorphism 
R\pi_{2,!}{\Bbb Q}_{\ell}
$$
in $D_{\rm c}^{\rm b}(Y,{\Bbb Q}_{\ell})$.  The composed isomorphisms 
$f_{1}\circ F_Y^{\ast} (f_{2})$ and $f_{2}\circ F_Y^{\ast} (f_{1})$ 
are the natural lifts of the Frobenius endomorphisms of $Y$ with 
respect to $k'$.  But, as the trace morphisms ${\rm tr}_{1}: 
R\pi_{1,!}{\Bbb Q}_{\ell}\longrightarrow {\Bbb Q}_{\ell}[-2r](-r)$ and 
${\rm tr}_{2}: R\pi_{2,!}{\Bbb Q}_{\ell}\longrightarrow {\Bbb 
Q}_{\ell}[-2r](-r)$ are both isomorphisms, we may view $f_{1}$ and 
$f_{2}$ as automorphisms of ${\Bbb Q}_{\ell}[-2r](-r)$.  It follows 
that $f_{1}$ and $f_{2}$ are induced by the multiplication on the 
constant sheaf ${\Bbb Q}_{\ell}$ by locally constant functions 
$\varphi_{1}: Y\rightarrow {\Bbb Q}_{\ell}$ and $\varphi_{2} 
:Y\rightarrow {\Bbb Q}_{\ell}$ such that the products 
$\varphi_{1}F_Y^{\ast}(\varphi_{2})$ and 
$\varphi_{2}F_Y^{\ast}(\varphi_{1})$ are both the constant function 
with value $q^{2r}$.

To conclude the proof of the proposition it is sufficient to check 
that the functions $\varphi_{1}$ and $\varphi_{2}$ are both constant 
with value $q^r$.  But this a direct consequence of the base change 
theorem in etale cohomology ([SGA 4] Exp. XVII, Thm.  5.2.6) and 
the next lemma.
\hfill\hfill$\square$

\th LEMMA 8.4
\enonce
Let $E_{1}$ and $E_{2}$ be two vector spaces of the same finite 
dimension $e$ over the finite field $k'$.  Let $F:E_{1}\rightarrow 
E_{2}$ be a $\sigma$-linear bijective map.  Let us view $E_{1}$ and 
$E_{2}$ as affine $k'$-schemes and $F$ as a finite $k'$-scheme 
morphism.  Then we have
$$
{\rm tr}_{1}\circ F^{\ast} =q^e{\rm tr}_{2}
$$
where ${\rm tr}_{i}:H_{\rm c}^{2e}(\overline k\otimes_{k'}E_{i},{\Bbb 
Q}_\ell )\buildrel{\sim}\over\longrightarrow {\Bbb Q}_\ell (-e)$ is 
the trace morphism and
$$
F^{\ast} : H_{\rm c}^{2e}(\overline k\otimes_{k'} E_{2},{\Bbb Q}_\ell 
)\rightarrow H_{\rm c}^{2e}(\overline k\otimes_{k'}E_{1},{\Bbb Q}_\ell 
)
$$
is induced by $F$.
\endth

\rem Proof
\endrem
The degree of the $k'$-scheme morphism $F: E_{1}\rightarrow E_{2}$ is 
$q^e$.
\hfill\hfill$\square$
\vskip 3mm

It is now clear that the Langlands-Shelstad conjecture, and more 
generally the equality
$$
\big|X^{F_{X}^{f}}\big|-\big|X'^{F_{X'}^{f}}\big|= 
q^{fr}\cdot\big|Y_{1}^{F_{Y_{1}}^{f}}\big|\cdot\big|Y_{2}^{F_{Y_{1}}^{f}}\big| 
$$
for every positive odd integer $f$, would be an immediate consequence 
of the Grothendieck fixed point formula (cf.  [Gr]) and the following 
conjecture:

\th CONJECTURE 8.5
\enonce
The two distinguished triangles
$$
R\Gamma (Y)[-2r](-r)\,\longmaprightover{(\pi_{1}^{!})^{-1}\circ 
j_{1,!}}{18mm}\,R\Gamma 
(X)\,\longmaprightover{i_{1}^{\ast}}{10mm}\,R\Gamma 
(X')\,\longmaprightover{\partial_{1}}{10mm}\,R\Gamma (Y)[1-2r](-r)
$$
and
$$
R\Gamma (Y)[-2r](-r)\,\longmaprightover{(\pi_{2}^{!})^{-1}\circ 
j_{2,!}}{18mm}\,R\Gamma 
(X)\,\longmaprightover{i_{2}^{\ast}}{10mm}\,R\Gamma 
(X')\,\longmaprightover{\partial_{2}}{10mm}\,R\Gamma (Y)[1-2r](-r)
$$
are identical.
\endth

To conclude this section we will discuss a possible homotopy argument 
for proving the equality of the restriction maps $i_{1}^{\ast}$ and 
$i_{2}^{\ast}$.
\vskip 3mm

The $k'$-schemes $X$, $X'$ and $Y_{i}$ are naturally embedded into 
Grassmann varieties ${\cal G}$, ${\cal G}'$ and ${\cal H}_{i}$.

More precisely, let ${\cal G}$ be the Grassmann variety of 
$(m_{1}+m_{2}+r)$-planes in the $2(m_{1}+m_{2}+r)$-dimensional 
$k'$-vector space $V_{1}\oplus V_{2}$ where we have set
$$
V_{i}=\left\{\matrix{\varpi_{E_{i}}^{-m_{i}-r'}{\cal 
O}_{E_{i}'}/\varpi_{E_{i}}^{m_{i}+r'}{\cal O}_{E_{i}'}&\hbox{ if 
}r=2r'\hbox{ is even,}\hfill\cr
\noalign{\medskip}
\varpi_{E_{i}}^{-m_{i}-r'-1}{\cal 
O}_{E_{i}'}/\varpi_{E_{i}}^{m_{i}+r'}{\cal O}_{E_{i}'}&\hbox{ if 
}r=2r'+1\hbox{ is odd,}\hfill\cr}\right.
$$
let ${\cal G}'$ be the Grassmann variety of $(m_{1}+m_{2}+r-1)$-planes 
in the $2(m_{1}+m_{2}+r-1)$-dimensional $k'$-vector space 
$V_{1}'\oplus V_{2}' $ where we have set
$$
V_{i}'=\left\{\matrix{\varpi_{E_{i}}^{-m_{i}-r'}{\cal 
O}_{E_{i}'}/\varpi_{E_{i}}^{m_{i}+r'-1}{\cal O}_{E_{i}'}&\hbox{ if 
}r=2r'\hbox{ is even,}\hfill\cr
\noalign{\medskip}
\varpi_{E_{i}}^{-m_{i}-r'}{\cal 
O}_{E_{i}'}/\varpi_{E_{i}}^{m_{i}+r'}{\cal O}_{E_{i}'}&\hbox{ if 
}r=2r'+1\hbox{ is odd,}\hfill\cr}\right.
$$
and let ${\cal H}_{i}$ be the Grassmann variety of $m_{i}$-planes in 
the $2m_{i}$-dimensional $k'$-vector space
$$
W_{i}=\varpi_{E_{i}}^{-m_{i}}{\cal 
O}_{E_{i}'}/\varpi_{E_{i}}^{m_{i}}{\cal O}_{E_{i}'}.
$$
The multiplications by $\varpi_{F}$ and $\gamma_{i}$ induce an 
endomorphism $\nu_{i}$ and an automorphism $u_{i}$ of $V_{i}$ (resp.  
$V_{i}'$, resp.  $W_{i}$). Then, if we set
$$
\nu =\nu_{1}\oplus\nu_{2}\quad \hbox{ and }\quad u=u_{1}\oplus u_{2}
$$
we have the obvious identifications:
$$
X=\{L\in {\cal G}\mid \nu (L)\subset L,~u(L)=L\hbox{ and }b_{1}\geq 
m_{1},b_{2}\geq m_{2}\}
$$
where $b_{1}={\rm dim}(L\cap (V_{1}\oplus (0)))$ and $b_{2}={\rm 
dim}(L\cap ((0)\oplus V_{2}))$,
$$
X'=\{L'\in {\cal G}'\mid \nu (L')\subset L',~u(L')=L'\hbox{ and } 
b_{1}'\geq m_{1},b_{2}'\geq m_{2}\}
$$
where $b_{1}'={\rm dim}(L'\cap (V_{1}'\oplus (0)))$ and $b_{2}'={\rm 
dim}(L'\cap ((0)\oplus V_{2}'))$, and
$$
Y_{i}=\{M_{i}\in {\cal H}_{i}\mid \nu_{i}(M_{i})\subset M_{i}\hbox{ 
and }u_{i}(M_{i})=M_{i}\}.
$$

\th CONJECTURE 8.6
\enonce
The restriction maps
$$
R\Gamma ({\cal G})\rightarrow R\Gamma (X),
$$
$$
R\Gamma ({\cal G}')\rightarrow R\Gamma (X')
$$
and
$$
R\Gamma ({\cal H}_{i})\rightarrow R\Gamma (Y_{i})
$$
induce epimorphisms on the cohomology.
\endth

\rem Remark $8.7$
\endrem
It follows from Conjecture $8.6$ that the cohomology complexes 
$R\Gamma (X)$, $R\Gamma (X)$ and $R\Gamma (X')$ should have all their 
odd cohomology groups equal to $(0)$ and that, for every even integer 
$n$ all the eigenvalues of Frobenius acting on the $n$-th cohomology 
group should be equal to $q^{{n\over 2}}$. 

In particular the two boundary maps $\partial_{1},\partial_{2}:R\Gamma 
(X')\rightarrow R\Gamma (Y)[1-2r](-r)$ should be zero and thus equal.
\vskip 3mm

For each flag
$$
F=((0)\subset F^{2}\subset F^{1}\subset V_{1}\oplus V_{2})
$$
of vector subspaces of $V_{1}\oplus V_{2}$ with $F^{1}$ of codimension 
$1$ and $F^{2}$ of dimension $1$ and for each isomorphism $\iota$ from 
the vector space $V_{1}'\oplus V_{2}'$ onto the vector space 
$F^{1}/F^{2}$, we have an obvious closed embedding
$$
i_{F,\iota}: {\cal G}'\hookrightarrow {\cal G}.
$$
As the $k'$-scheme of pairs $(F,\iota )$ is connected, the restriction maps
$$
i_{F,\iota}^{\ast}: R\Gamma ({\cal G})\rightarrow R\Gamma ({\cal G}')
$$
are all equal. We will simply denote by $i^\ast$ the common value of these
restriction maps. 

Let us denote by $N_{i}$ the regular nilpotent endomorphism of 
$V_{i}$ which is induced by the multiplication by $\varpi_{E_{i}}$. By 
definition we have
$$
V_{i}'=\left\{\matrix{{\rm Coker}(N_{i}) 
\,\longtwoheadleftarrowover{}{6mm}\, V_{i}&\hbox{if }r\hbox{ is 
even}\hfill\cr
\noalign{\medskip}
{\rm Im}(N_{i})\,\longhookrightarrowover{}{6mm}\, V_{i}&\hbox{if }r\hbox{ 
is odd}\hfill\cr}\right.
$$
Moreover the endomorphism $N_{i}$ of $V_{i}$ induces an isomorphism of 
${\rm Coker}(N_{i})$ onto ${\rm Im}(N_{i})$ and we may identify these 
two subquotients of $V_{i}$.  If
$$
F=((0)\subset {\rm Ker}(N_{1})\oplus (0)\subset V_{1}\oplus {\rm 
Im}(N_{2}) \subset V_{1}\oplus V_{2})
$$
and if $\iota$ is the isomorphism
$$
V_{1}'\oplus V_{2}'\isomorphism (V_{1}\oplus {\rm Im}(N_{2}))/({\rm 
Ker}(N_{1})\oplus (0))
$$
which is induced by the above identifications, the closed embedding
$i_{F,\iota}$ maps $X'$ into $X$ and extends the closed embedding
$i_{1}:X'\hookrightarrow X$. Similarly, if
$$
F=((0)\subset (0)\oplus {\rm Ker}(N_{2})\subset {\rm Im}(N_{1})\oplus V_{2}
\subset V_{1}\oplus V_{2})
$$
and $\iota$ is the isomorphism
$$
V_{1}'\oplus V_{2}'\isomorphism ({\rm Im}(N_{1})\oplus 
V_{2})/((0)\oplus {\rm Ker}(N_{2}))
$$
which is induced by the above identifications, the closed embedding
$i_{F,\iota}$ maps $X'$ into $X$ and extends the closed embedding
$i_{2}:X'\hookrightarrow X$. In particular, for $i=1,2$, we have a
commutative diagram
$$\diagram{
R\Gamma ({\cal G})& \maprightover{i^\ast}&R\Gamma ({\cal G}')\cr
\mapdownleft{}&&\mapdownright{}\cr
R\Gamma (X)& \maprightunder{i_i^\ast} &R\Gamma (X')\cr}
$$
and the equality $i_{1}^{\ast}=i_{2}^\ast$ follows from Conjecture $8.6$.

\vskip 5mm

\centerline{\bf 9. $U(1,1)$}
\vskip 5mm

In this section we assume that $n_{1}=n_{2}=1$.  Replacing 
$(\gamma_{1},\gamma_{2})$ by $(1,\gamma_{1}^{-1}\gamma_{2})$ we may 
also assume that $\gamma_{1}=1$.  Then $\gamma_{2}-1$ is of valuation 
$r$ in $F'$.

Let $K$ be a field extension of $k'$.  If $L^{\pm}$ is a 
$(K\otimes_{k'}{\cal O}_{F'})$-lattice in $F'\oplus F'$ of index 
$\delta^{\pm}$ ($\delta^{+}=0$ and $\delta^{-}=1$) we necessarily have
$$
B_{i}^{\pm}=\varpi_{F}^{-b_{i}^{\pm}}(K\otimes_{k'}{\cal 
O}_{F'})\subset 
C_{i}^{\pm}=\varpi_{F}^{-c_{i}^{\pm}}(K\otimes_{k'}{\cal 
O}_{F'})\subset K\otimes_{k'}F'
$$
for some integers $b_{i}^{\pm}\leq c_{i}^{\pm}$ such that
$$
b_{1}^{\pm}+c_{2}^{\pm}= b_{2}^{\pm}+c_{1}^{\pm}=\delta^{\pm}.
$$
Therefore the conditions
$$
(1\otimes (\gamma_{1}\oplus\gamma_{2}))L^{\pm}=L^{\pm}
$$
and
$$
c_{1}^{\pm}-b_{1}^{\pm}=c_{2}^{\pm}-b_{2}^{\pm}\leq r
$$
are equivalent.

If $r=2r'$ is even it follows that
$$
X=X^{+}[-r',-r']
$$
is simply the $k'$-scheme of ${\cal O}_{F'}$-lattices $L^{+}$ in 
$F'\oplus F'$ satisfying the conditions
$$\left\{\matrix{
\varpi_{F}^{r'}{\cal O}_{F'}\oplus\varpi_{F}^{r'}{\cal O}_{F'}\subset 
L^{+}\subset\varpi_{F}^{-r'}{\cal O}_{F'}\oplus\varpi_{F}^{-r'}{\cal 
0}_{F'}\hfill\cr
\noalign{\medskip}
{\rm ind}(L^{+})=0,\hfill\cr}\right.
$$
and that
$$
X'=X^{-}[1-r',1-r']
$$
is simply the $k'$-scheme of ${\cal O}_{F'}$-lattices $L^{-}$ in 
$F'\oplus F'$ satisfying the conditions
$$\left\{\matrix{
\varpi_{F}^{r'-1}{\cal O}_{F'}\oplus\varpi_{F}^{r'-1}{\cal O}_{F'}\subset 
L^{-}\subset\varpi_{F}^{-r'}{\cal O}_{F'}\oplus\varpi_{F}^{-r'}{\cal 
0}_{F'}\hfill\cr
\noalign{\medskip}
{\rm ind}(L^{-})=1.\hfill\cr}\right.
$$

Similarly, if $r=2r'+1$ is odd it follows that
$$
X=X^{-}[-r',-r']
$$
is simply the $k'$-scheme of ${\cal O}_{F'}$-lattices $L^{-}$ in 
$F'\oplus F'$ satisfying the conditions
$$\left\{\matrix{
\varpi_{F}^{r'}{\cal O}_{F'}\oplus\varpi_{F}^{r'}{\cal O}_{F'}\subset 
L^{+}\subset\varpi_{F}^{-r'-1}{\cal O}_{F'}\oplus\varpi_{F}^{-r'-1}{\cal 
0}_{F'}\hfill\cr
\noalign{\medskip}
{\rm ind}(L^{-})=1,\hfill\cr}\right.
$$
and that
$$
X'=X^{+}[-r',-r']
$$
is simply the $k'$-scheme of ${\cal O}_{F'}$-lattices $L^{+}$ in 
$F'\oplus F'$ satisfying the conditions
$$\left\{\matrix{
\varpi_{F}^{r'}{\cal O}_{F'}\oplus\varpi_{F}^{r'}{\cal O}_{F'}\subset 
L^{+}\subset\varpi_{F}^{-r'}{\cal O}_{F'}\oplus\varpi_{F}^{-r'}{\cal 
0}_{F'}\hfill\cr
\noalign{\medskip}
{\rm ind}(L^{+})=0.\hfill\cr}\right.
$$

If $r=2r'$ is even (resp.  $r=2r'+1$ is odd) let us denote by $V$ 
the $r$-dimensional $k'$-vector space
$$
\varpi_{F}^{-r'}{\cal O}_{F'}/\varpi_{F}^{r'}{\cal O}_{F'}
$$
(resp.
$$
\varpi_{F}^{-r'-1}{\cal O}_{F'}/\varpi_{F}^{r'}{\cal O}_{F'}\,)
$$
and by $N$ the nilpotent endomorphism of $V$ which is induced by the 
multiplication by $\varpi_{F}$.  Then $X$ may be identified with the 
$k'$-scheme of $k'$-vector subspaces
$$
A\subset V\oplus V
$$
satisfying the conditions
$$\left\{\matrix{
(N\oplus N)(A)\subset A\hfill\cr
\noalign{\medskip}
{\rm dim}(A)=r.\hfill\cr}\right.
$$
Similarly, if $r=2r'$ is even (resp.  $r=2r'+1$ is odd) let us denote 
by $V'$ the $(r-1)$-dimensional $k'$-vector space
$$
\varpi_{F}^{-r'}{\cal O}_{F'}/\varpi_{F}^{r'-1}{\cal O}_{F'}
$$
(resp.
$$
\varpi_{F}^{-r'}{\cal O}_{F'}/\varpi_{F}^{r'}{\cal O}_{F'}\,)
$$
and by $N'$ the nilpotent endomorphism of $V'$ which is induced by the 
multiplication by $\varpi_{F}$.  Then $X'$ may be identified with the 
$k'$-scheme of $k'$-vector subspaces
$$
A'\subset V'\oplus V'
$$
satisfying the conditions
$$\left\{\matrix{
(N'\oplus N')(A')\subset A'\hfill\cr
\noalign{\medskip}
{\rm dim}(A)=r-1.\hfill\cr}\right.
$$

Moreover under these identifications the closed embeddings $i_{1}$ and 
$i_{2}$ defined in (7.1) may be described in the following way.  Let us 
consider the diagram
$$\diagram{
&&V\oplus V&&\cr
\noalign{\vskip 5mm}
&&&&\cr
\noalign{\vskip 5mm}
V\oplus V'&&&&V'\oplus V\cr
\noalign{\vskip 5mm}
&&&&\cr
\noalign{\vskip 5mm}
&&V'\oplus V'&&\cr
\put (7,29){\scriptstyle 1\oplus e}
\arrow(10,25)\dir(1,1)\length{7}
\put (33.7,29){\scriptstyle e\oplus 1}
\arrow(34.7,25)\dir(-1,1)\length{7}
\put (7,10){\scriptstyle p\oplus 1}
\arrow(10,15)\dir(1,-1)\length{7}
\put (33.7,10){\scriptstyle 1\oplus p}
\arrow(34.7,15)\dir(-1,-1)\length{7}}
$$
where $e:V'\,\longhookrightarrowover{}{6mm}\,V$ is the embedding
$$
\varpi_{F}^{-r'}{\cal O}_{F'}/\varpi_{F}^{r'-1}{\cal O}_{F'} 
\,\longhookrightarrowover{}{6mm}\,\varpi_{F}^{-r'}{\cal 
O}_{F'}/\varpi_{F}^{r'}{\cal O}_{F'}
$$
which is induced by the multiplication by $\varpi_{F}$ (resp.  the 
canonical embedding
$$
\varpi_{F}^{-r'}{\cal O}_{F'}/\varpi_{F}^{r'}{\cal O}_{F'} 
\,\longhookrightarrowover{}{6mm}\,\varpi_{F}^{-r'-1}{\cal 
O}_{F'}/\varpi_{F}^{r'}{\cal O}_{F'}\,)
$$
and where $p:V\,\longtwoheadrightarrowover{}{6mm}\,V'$ is the 
canonical projection
$$
\varpi_{F}^{-r'}{\cal O}_{F'}/\varpi_{F}^{r'}{\cal O}_{F'} 
\,\longtwoheadrightarrowover{}{6mm}\,\varpi_{F}^{-r'}{\cal 
O}_{F'}/\varpi_{F}^{r'-1}{\cal O}_{F'}
$$
(resp.  the projection
$$
\varpi_{F}^{-r'}{\cal O}_{F'}/\varpi_{F}^{r'}{\cal O}_{F'} 
\,\longtwoheadrightarrowover{}{6mm}\,\varpi_{F}^{-r'}{\cal 
O}_{F'}/\varpi_{F}^{r'-1}{\cal O}_{F'}
$$
which is induced by the multiplication by $\varpi_{F}$).  Let us 
remark that $e$ identifies $V'$ with ${\rm Im}(N)$ and $N'$ with the 
restriction of $N$ to its image, and that $p$ identifies $V'$ with 
${\rm Coker}(N)$ and $N'$ with the regular nilpotent endomorphism 
induced by $N$ on its cokernel.  Then
$$
i_{1}(A')=(1\oplus e)((p\oplus 1)^{-1}(A'))
$$
and
$$
i_{2}(A')=(e\oplus 1)((1\oplus p)^{-1}(A')).
$$

\th THEOREM 9.1
\enonce
The two restriction maps
$$
i_{1}^{\ast} ,i_{2}^{\ast} : R\Gamma (X)\rightarrow R\Gamma (X')
$$
are equal.
\endth

\rem Proof
\endrem
As we have explained at the end of Section 8 it is sufficient to 
prove that the restriction maps
$$
R\Gamma ({\cal G})\rightarrow R\Gamma (X)
$$
where ${\cal G}$ is the Grassmann variety of $r$-planes in $V\oplus 
V$, and
$$
R\Gamma ({\cal G}')\rightarrow R\Gamma (X')
$$
where ${\cal G}'$ is the Grassmann variety of $(r-1)$-planes in 
$V'\oplus V'$, induce epimorphisms on the cohomology.  But this has 
been proved by Hotta and Shimomura as a consequence of a result of 
Spaltenstein (cf.  [Ho-Sh] Lemma 8.1 and its proof).
\hfill\hfill$\square$
\vskip 3mm

\th COROLLARY 9.2
\enonce
Conjecture $8.5$ holds in the case $n_{1}=n_{2}=1$.  In particular, 
for any finite extension $k_f$ of degree $f$ of $k$ we have the 
identity
$$
|{\cal X}^{+} (k_f)|-|{\cal X}^{-}(k_f)|=(-1)^rq^{rf}\cdot |{\cal Y}(k_f)|
$$
with $|{\cal Y}(k_f)|=1$.
\endth

\rem Proof
\endrem
We know already that the two restriction maps $i_{1}^{\ast}$ and 
$i_{2}^{\ast}$ are equal.

As $Y$ is clearly equal to ${\rm Spec}(k')$ ($m_{1}=m_{2}=0$) it 
follows from Section 7 that $U_{i}=X-i_i(X')$ is a standard affine 
space of dimension $r$ over $k'$.  But $X'$ is isomorphic to the 
variety $X$ after having replaced $r$ by $r-1$.  Therefore an obvious 
induction on $r$ shows that $X$ is a disjoint union of standard affine 
spaces, one for each dimension between $0$ and $r$.  Consequently 
either the source or the target of $\partial_{i}$ is zero and we get
$$
\partial_{1}=\partial_{2}=0.
$$

Finally let us consider the two Gysin maps $j_{1,!}$ and $j_{2,!}$.  
The only degree where these two maps can be non zero is the top degree 
$2r$.  In this degree the $\ell$-adic cohomology groups with compact 
supports of $U_{1}$, $U_{2}$ and $X$ can be all canonically identified 
with ${\Bbb Q}_\ell (-r)$.  Under these identifications 
$(\pi_{1}^{!})\circ j_{1,!}$ and $(\pi_{2}^{!})\circ j_{2,!}$ become 
the identity of ${\Bbb Q}_\ell (-r)$ and are thus equal.
\hfill\hfill$\square$
\vskip 3mm

In the $U(1,1)$ case the closed embeddings $i_{1},i_{2}:X'\, 
\longhookrightarrowover{}{6mm}\,X$ are homotopic.  We will conclude 
this section by constructing such an homotopy and giving another proof 
of Theorem 9.1.  We introduce the projective $k'$-scheme 
$\widetilde{\widetilde{X}}$ of partial flags
$$
L\subset A\subset H\subset V\oplus V
$$
of vector subspaces such that $L$, $A$ and $H$ are of dimensions $1$, 
$r$ and $2r-1$ respectively, and are stable by $N \oplus N$.  
Obviously, for each point $L\subset A\subset H$ in
$\widetilde{\widetilde{X}}$ we have
$$
L\subset {\rm Ker}(N)\oplus {\rm Ker}(N)
$$
and
$$
H\supset {\rm Im}(N)\oplus {\rm Im}(N).
$$
By forgetting either $A$ or $(L,H)$ we get two $k'$-scheme morphisms
$$\diagram{
\widetilde{\widetilde{X}} &\kern 
-6mm\longmaprightover{\tilde{g}}{10mm}\kern -2mm&X\cr
\mapdownleft{\tilde{f}}&&\cr
S\times S'&&\cr}\leqno (9.3)
$$
where $S$ is the projective line of lines $L$ in the two dimensional 
$k'$-vector space ${\rm Ker}(N)\oplus {\rm Ker}(N)$ and $S'$ is the 
projective line of hyperplanes $H$ in $V\oplus V$ which contain the 
codimension $2$ subspace ${\rm Im}(N)\oplus {\rm Im}(N)$. We have a 
closed embedding
$$
S\hookrightarrow S\times S',~L \mapsto (L,(N\oplus N )^{-(r-1)}(L)).
$$
Let us denote by
$$\diagram{
\widetilde{X}&\kern -2mm\longmaprightover{g}{10mm}\kern -2mm&X\cr
\mapdownleft{f}&&\cr
S&&\cr}\leqno (9.4)
$$
the inverse image of the diagram $(9.3)$ by that closed embedding. 

In $S$ we have two $k'$-rational points $s_{1}$ and $s_{2}$ given by 
the lines ${\rm Ker}(N)\oplus (0)$ and $(0)\oplus {\rm Ker}(N)$ in 
${\rm Ker}(N)\oplus {\rm Ker}(N)$.  The restriction of $g$ to the 
fiber of $f$ over $s_i$ is an isomorphism onto $i_i(X')$.  We will 
denote by $\tilde{i}_i:X'\hookrightarrow \widetilde{X}$ the obvious 
lifting of $i_i$ to $\widetilde{X}$ with image that fiber.  More 
generally, for any point $s$ in $S$, corresponding to a line $L$ in 
${\rm Ker}(N)\oplus {\rm Ker}(N)$, the restriction of $g$ to the fiber 
of $f$ over $s$ is an isomorphism onto the closed subset
$$
\{A\in X\mid L\subset A\subset (N\oplus N)^{-(r-1)}(L)\}
$$
of $X$.

\th PROPOSITION 9.5
\enonce
{\rm (i)} The projective morphism $g$ is birational.  More precisely 
it is an isomorphism over the complementary subset in $X$ of the 
closed subset
$$
i_{1}(X')\cap i_{2}(X')=\{A\in X\mid {\rm Ker}(N)\oplus {\rm 
Ker}(N)\subset A\subset {\rm Im}(N)\oplus {\rm Im}(N)\}
$$
and its restriction to this closed subset is a trivial fibration in 
projective lines.

\decale{\rm (ii)} The projective morphism $f$ is a locally trivial 
fibration for the Zariski topology, with all its fibers isomorphic to 
$X'$.  More precisely, if $\{i,j\}=\{1,2\}$ the restriction of $f$ to 
the Zariski open subset $S-\{s_i\}\subset S$ is isomorphic to the 
canonical projection $(S-\{s_i\})\times X'\rightarrow S-\{s_i\}$ by an 
isomorphism which exchanges the closed embeddings
$$
\tilde{i}_j:X'\hookrightarrow f^{-1}(S-\{s_i\})
$$
and
$$
X'\cong \{s_j\}\times X'\subset (S-\{s_i\})\times X'.
$$
Moreover the gluing datum between the above trivializations of $f$ 
over the open subsets $S-\{s_{1}\}$ and $S-\{s_{2}\}$ is given by a 
morphism
$$
S-\{s_{1},s_{2}\}\rightarrow {\rm Aut}(X')
$$
which factors through the canonical morphism from the centralizer of 
$N'\oplus N'$ in ${\rm GL}(V'\oplus V')$ to the automorphism group of 
$X'$.
\endth

\rem Proof
\endrem
Part (i) is obvious.

Let us prove Part (ii).  The open subset $S-\{s_{1}\}\subset S$ may be 
identified with the affine line of linear maps $\varphi: {\rm 
Ker}(N)\rightarrow {\rm Ker}(N)$.  For example $s_{2}$ corresponds to 
$\varphi=0$.  The fiber of $f$ over a point $\varphi$ is the closed 
subvariety
$$
X(\varphi )=\{A\in X\mid L(\varphi )\subset A\subset H(\varphi )\}
$$
of $X$ where
$$
L(\varphi )=\{(\varphi(x_{2}),x_{2})\mid x_{2}\in {\rm Ker}(N)\}
$$
and
$$
H(\varphi )=\{(x_{1},x_{2})\in V\oplus V\mid N^{r-1}(x_{1})= 
\varphi\circ N^{r-1}(x_{2})\}.
$$

Let us fix a cyclic vector $v\in V$ for $N$.  Then $N^{r-1}(v)$ is a 
basis of ${\rm Ker}(N)$ and any linear map $\varphi$ as above maps 
$N^{r-1}(v)$ onto $\lambda (\varphi )N^{r-1}(v)$ for a unique scalar 
$\lambda (\varphi )$.  Let us denote by $\widetilde{\varphi}$ the 
unique endomorphism of $V$ such that
$$
\widetilde{\varphi}(v)=\lambda (\varphi )v
$$
and
$$
N\circ\widetilde{\varphi}=\widetilde{\varphi}\circ N.
$$
Obviously $\widetilde{\varphi}$ extends $\varphi$ and we may identify 
$H(\varphi )/L(\varphi )$ with ${\rm Im}(N)\oplus {\rm Im}(N)=V'\oplus 
V'$ by sending $(x_{1},x_{2})\in H(\varphi )$ onto
$$
(x_{1}-\widetilde{\varphi}(x_{2}),N(x_{2})).
$$
As this identification exchanges the endomorphisms which are induced 
by $N\oplus N$ it induces an isomorphism from $X(\varphi )$ onto $X'$.  
Letting $\varphi$ vary we get an isomorphism from the restriction of 
$f$ to $S-\{s_{1}\}$ onto the canonical projection 
$(S-\{s_{1}\})\times X'\rightarrow S-\{s_{1}\}$.  Clearly this 
isomorphism exchanges the closed embeddings $\tilde{i}_{2} 
:X'\hookrightarrow f^{-1}(S-\{s_{1}\})$ and $X'\cong \{s_{2}\}\times 
X'\subset (S-\{s_{1}\})\times X'$.

Similarly we identify $S-\{s_{2}\}$ with the affine line of linear 
maps $\psi: {\rm Ker}(N)\rightarrow {\rm Ker}(N)$ and we construct an 
isomorphism from the restriction of $f$ to $S-\{s_{2}\}$ onto the 
canonical projection $(S-\{s_{2}\})\times X'\rightarrow S-\{s_{2}\}$.  
With obvious notations it is given by sending $(x_{1},x_{2})\in 
H(\psi )$ 
onto
$$
(N(x_{1}),x_{2}-\widetilde{\psi}(x_{1})).
$$

Over $S-\{s_{1},s_{2}\}$ the gluing datum between the two 
trivializations of $f$ is given by the morphism
$$
S-\{s_{1},s_{2}\}\rightarrow {\rm Aut}(X')
$$
which sends $\varphi =\psi^{-1}$ onto the automorphism of $X'$ which 
is induced by the automorphism
$$
(x_{1},x_{2})\mapsto 
(-\widetilde{\varphi}(x_{2}),\widetilde{\varphi}^{-1}(x_{1}) 
+N'(x_{2}))
$$
of $V\oplus V$.
\hfill\hfill$\square$
\vskip 3mm

\rem Another proof of Theorem $9.1$
\endrem
It is sufficient to check that the two restriction maps
$$
\tilde{i}_{1}^{\ast} ,\tilde{i}_{2}^{\ast} : R\Gamma 
(\widetilde{X})\rightarrow R\Gamma (X')
$$
are equal (we have $i_{i}^{\ast}=\tilde{i}_{i}^{\ast}\circ g^{\ast}$).  
By the Leray spectral sequence we have
$$
R\Gamma (\widetilde{X}) =R\Gamma (\overline k\otimes_{k'}S, Rf_\ast 
{\Bbb Q}_\ell ).
$$
But the proposition implies that

\decale{\rm (i)} the complex $Rf_\ast {\Bbb Q}_\ell$ is constant with 
value $R\Gamma (X')$ on the open subsets $S-\{s_{1}\}$ and 
$S-\{s_{2}\}$ of $S$ and we can choose the trivializations in such way 
that $\tilde{i}_{1}^{\ast}$ and $\tilde{i}_{2}^{\ast}$ are both 
induced by the identity of $R\Gamma (X')$,

\decale{\rm (ii)} the gluing datum for $Rf_\ast {\Bbb Q}_\ell$ on 
$S-\{s_{1},s_{2}\}$ is induced by the identity of $R\Gamma (X')$ as the 
centralizer of $N'\oplus N'$ in ${\rm GL}(V'\oplus V')$ is connected.
\vskip 1mm

The theorem follows.
\hfill\hfill$\square$
\vskip 30mm

\centerline{\bf 10. Remarks and examples}
\vskip 5mm

(10.1) In general the conductors $m_{1}$, $m_{2}$ and the order of the 
resultant $r$ are difficult to compute.  Nevertheless there are some 
interesting particular cases where they admit simple formulae.

In general we have
$$
m_{i}=v_{E_{i}'}(P_{i}'(\gamma_{i}))-\delta (E_{i}'/F')
$$
and
$$
r=v_{E_{i}'}(P_{j}(\gamma_{i})).
$$
Here $v_{E_{i}'}:E_{i}'\rightarrow {\Bbb Z}\cup\{+\infty\}$ is the
discrete valuation of $E_{i}'$, $P_{i}'(T)$ is the derivative of the 
minimal polynomial $P_{i}(T)$ of $\gamma_{i}$ over $F'$, $\delta 
(E_{i}'/F')$ is the order of the different of the totally ramified 
extension $E_{i}'/F'$ and $\{i,j\}=\{1,2\}$ (cf.  [Se] Ch.  III, \S 6, 
Corollaire 1).  Moreover, if $n_{i}$ is prime to the characteristic of 
$k$ we have
$$
\delta (E_{i}'/F')=n_{i}-1
$$
(cf. [Se] Ch.  III, \S 6, Proposition 13).

Therefore, if we denote by $a_{i}$ the constant term of the expansion of 
$\gamma_{i}$ as power series in $\varpi_{E_{i}}$ with coefficients in 
$k'$ and we put
$$
v_{i}=v_{E_{i}'}(\gamma_{i}-a_{i})
$$
we have:

\decale{\rm (i)} $m_{i}=0$ and $r=v_{E_{j}'}(\gamma_{j}-\gamma_{i})$ if 
$n_{i}=1$,

\decale{\rm (ii)} $m_{i}=(n_{i}-1)(v_{i}-1)$ if $n_{i}>1$ is prime to 
the characteristic of $k$ and $v_{i}$ is prime to $n_{i}$,

\decale{\rm (ii)} $r\geq {\rm Inf}(n_{1}v_{2},n_{2}v_{1})$ with equality 
if $n_{1}v_{2}\not=n_{2}v_{1}$.
\vskip 5mm

(10.2) {\it The case $r=0$}.  In this case the $k'$-schemes $X^{+}$ 
and $X^{-}$ are not connected: they are disjoint sums
$$
X^{+}=\coprod_{n\in {\Bbb Z}}X^{+}[-n,n]\quad\hbox{ and 
}\quad X^{-}=\coprod_{n\in {\Bbb Z}}X^{-}[1-n,n]
$$
where each component $X^{+}[-n,n]$ or $X^{-}[1-n,n]$ is isomorphic to 
$Y=Y_{1}\times_{k'}Y_{2}$.  We have $X=X^{+}[0,0]=Y$ with 
$F_{X}=F_{Y}$ and $X'=X^{-}[1,1]=\emptyset$, and Conjecture 8.5 and 
the Langlands-Shelstad conjecture are obvious.
\vskip 5mm

(10.3) {\it The case $n_{1}=1<n_{2}$ and $r>0$}.  In this case, we 
have $m_{1}=0$ and $Y_{1}={\rm Spec}(k')$, we may assume that 
$\gamma_{1}=1$ and thus $\gamma_{2}-1$ is of order $r$ in $E_{2}'$.

On $Y=Y_{2}$ we have the rank $r$ vector bundles
$$
\pi_{1}:U_{1}\rightarrow Y
$$
and
$$
\pi_{2}:U_{2}\rightarrow Y
$$
defined by
$$
\pi_{2}^{-1}(M_{2})=(M_{2}/(\gamma_{2}-1)M_{2})^{\vee}
$$
where $(M_{2}/(\gamma_{2}-1)M_{2})^{\vee}$ is the dual of the 
$r$-dimensional vector space $M_{2}/(\gamma_{2}-1)M_{2}$ and
$$
\pi_{2}^{-1}(M_{2})=(\gamma_{2}-1)^{-1}M_{2}/M_{2}
$$
respectively.  For $i=1,2$ and for every integer $j$ the rank 
$r$-vector bundle $\pi_{i,j}^{\pm}:U_{i,j}^{\pm}\rightarrow Y$ which 
has been defined in Section 5 is isomorphic to 
$\pi_{i}:U_{i}\rightarrow Y$, the isomorphism being given by
$$
(B_{1}^{\pm},C_{2}^{\pm},f_{1}^{\pm})\mapsto 
(M_{2}=\varpi_{E_{2}}^{c_{2}^{\pm}}C_{2}^{\pm},x_{1}^{\vee})
$$
where $x_{1}^{\vee}$ is induced by the image of the $k'$-linear form
$$
F'/B_{1}^{\pm}\rightarrow k',\quad\alpha\mapsto {\rm 
Res}(\alpha\varpi_{F}^{b_{1}^{\pm}-1}{\rm d}\varpi_{F})
$$
by the transpose of the morphism
$$
f_{1}^{\pm}\circ\varpi_{E_{2}}^{-c_{2}^{\pm}}: M_{2}\isomorphism 
C_{2}^{\pm}\,\longmaprightover{}{6 mm}\,F'/B_{1}^{\pm},
$$
and
$$
(B_{2}^{\pm},C_{1}^{\pm},f_{2}^{\pm})\mapsto 
(M_{2}=\varpi_{E_{2}}^{b_{2}^{\pm}}B_{2}^{\pm},x_{2})
$$
where $x_{2}$ is the image of the vector
$$
f_{2}^{\pm}(\varpi_{F}^{-c_{1}^{\pm}})\in E_{2}'/B_{2}^{\pm}
$$
by the isomorphism $\varpi_{E_{2}}^{b_{2}^{\pm}}: E_{2}'/B_{2}^{\pm} 
\isomorphism E_{2}'/M_{2}$.
\vskip 5mm

{\it In the next examples we assume that $n_{1}$ and $n_{2}$ are 
prime to the characteristic of $k$. We choose the uniformizers 
$\varpi_{F}$ and $\varpi_{E_{i}}$ in such way that 
$$
\varpi_{F}=\alpha_{i}\varpi_{E_{i}}^{n_{i}}
$$
in $E_{i}'$ for some $\alpha_{i}\in k'$.}
\vskip 5mm

(10.4) {\it The case $n_{1}=1$, $\gamma_{1}=1$, $n_{2}=3$, $r=2$ and 
${\rm Char}(k)>3$}.  In this case $\gamma_{2}-1$ is of order $2$ in 
$E_{2}'$ and $m_{2}=2$.  The $k'$-scheme $Y_{2}=Z_{2}$ is a projective 
line (cf. Remark 3.14).

Then the stratification of $X=X^{+}[-1,-1]$ by the values of $b_{1}^{+}$ and 
$b_{2}^{+}$ can be symbolically represented by the following diagram
$$
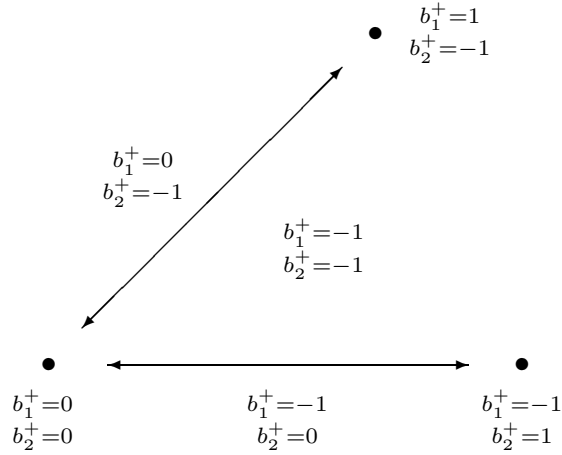
\diagram{
\noalign{\vskip 40mm}
\put(-26.4,-0.75){\bullet}
\put(26,-0.75){\bullet}
\put(9.8,36){\bullet}
\put(14,36){b_{1}^{+}=1\atop b_{2}^{+}=-1}
\put(-20,20){b_{1}^{+}=0\atop b_{2}^{+}=-1}
\put(0,12){b_{1}^{+}=-1\atop b_{2}^{+}=-1}
\put(-30,-7){b_{1}^{+}=0\atop b_{2}^{+}=0}
\put(-4,-7){b_{1}^{+}=-1\atop b_{2}^{+}=0}
\put(22,-7){b_{1}^{+}=-1\atop b_{2}^{+}=1}
\arrow(-8,19.15)\dir(-1,-1)\length{15}
\arrow(-8,18)\dir(1,1)\length{15}
\arrow(1,0)\dir(-1,0)\length{20}
\arrow(1,0)\dir(1,0)\length{20}
\noalign{\vskip 10mm}
}
$$
and the stratification of $X'=X^{-}[0,0]$ by the values of $b_{1}^{-}$ and 
$b_{2}^{-}$ can be symbolically represented by the following diagram
$$
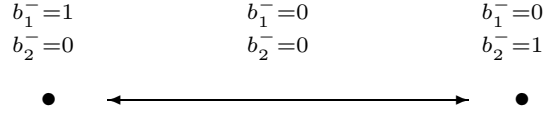
\diagram{
\noalign{\vskip 10mm}
\put(-30,7){b_{1}^{-}=1\atop b_{2}^{-}=0}
\put(-4,7){b_{1}^{-}=0\atop b_{2}^{-}=0}
\put(22,7){b_{1}^{-}=0\atop b_{2}^{-}=1}
\put(-26.4,-0.75){\bullet}
\put(26,-0.75){\bullet}
\arrow(1,0)\dir(-1,0)\length{20}
\arrow(1,0)\dir(1,0)\length{20}
\noalign{\vskip 5mm}
}
$$
In these representations the closed embeddings $i_{1}$ (resp. $i_{2}$) 
are given by
$$
(b_{1}^{-},b_{2}^{-})\mapsto 
(b_{1}^{+},b_{2}^{+})=(b_{1}^{-},b_{2}^{-}-1)
$$
(resp.
$$
(b_{1}^{-},b_{2}^{-})\mapsto 
(b_{1}^{+},b_{2}^{+})=(b_{1}^{-}-1,b_{2}^{-})\,).
$$

Let $V_{2}$ be the $6$-dimensional $k'$-vector space 
$$
V_{2}=\varpi_{E_{2}}^{-3}{\cal O}_{E_{i}'}/\varpi_{E_{2}}^{3}{\cal 
O}_{E_{i}'}
$$
and let $N_{2}$ and $\nu_{2}$ be the nilpotent endomorphisms of 
$V_{2}$ which are induced by the multiplication by $\varpi_{E_{2}}$ 
and $\varpi_{F}$, and let $u_{2}$ be the automorphism of $V_{2}$ which 
is induced by the multiplication by $\gamma_{2}$.

The stratum $(b_{1}^{+}=-1,b_{2}^{+}=-1)$ is isomorphic to the 
$k'$-scheme of triples
$$
(B_{2}\subset C_{2},x_{2})
$$
where $B_{2}$ and $C_{2}$ are vector subspaces of $V_{2}$ of dimension 
$2$ and $4$ respectively such that
$$
{\rm Ker}(N_{2})\subset B_{2}\subset {\rm Ker}(N_{2}^{3})\subset 
C_{2}\subset {\rm Ker}(N_{2}^{5})
$$
and
$$
(\gamma_{2}-1)(C_{2})=B_{2},
$$
and where $x_{2}$ is a vector in $C_{2}/B_{2}$ such that 
$$
\nu_{2}(x_{2})\not=0.
$$
Obviously the $k'$-scheme of pairs $(B_{2}\subset C_{2})$ as above is a 
projective line over $k'$ with some marked point $\infty :=({\rm 
Ker}(N_{2}^{2})\subset {\rm Ker}(N_{2}^{4}))$. For a given pair 
$(B_{2}\subset C_{2})$ there exist vectors $x_{2}\in C_{2}/B_{2}$ such that 
$\nu_{2}(x_{2})\not=0$ if and only if $(B_{2}\subset C_{2})$ is not equal 
to $\infty$.  Moreover, if this holds, the scheme of $x_{2}\in 
C_{2}/B_{2}$ such that $\nu_{2}(x_{2})\not=0$ is isomorphic to ${\Bbb 
A}^{1}\times {\Bbb G}_{{\rm m}}$.

The stratum $(b_{1}^{+}=0,b_{2}^{+}=-1)$ is isomorphic to the 
$k'$-scheme of triples
$$
(B_{2}\subset C_{2},x_{2})
$$
where $B_{2}$ and $C_{2}$ are vector subspaces of $V_{2}$ of dimension 
$2$ and $3$ respectively, such that
$$
{\rm Ker}(N_{2})\subset B_{2}\subset {\rm Ker}(N_{2}^{3})
$$
and
$$
{\rm Ker}(N_{2}^{2})\subset C_{2}\subset {\rm Ker}(N_{2}^{4})
$$
(the condition
$$
(\gamma_{2}-1)(C_{2})\subset B_{2}
$$
is automatic), and where $x_{2}$ is a non zero vector in 
$C_{2}/B_{2}$. If $B_{2}\not= {\rm Ker}(N_{2}^{2})$ we necessarily 
have
$$
C_{2}=B_{2}+{\rm Ker}(N_{2}^{2})={\rm Ker}(N_{2}^{3})
$$
and if $B_{2}={\rm Ker}(N_{2}^{2})$, $C_{2}$ may vary freely in the 
projective line
$$
{\rm Ker}(N_{2}^{2})\subset C_{2}\subset {\rm Ker}(N_{2}^{4}).
$$
Similarly, if $C_{2}\not= {\rm Ker}(N_{2}^{3})$ we necessarily 
have
$$
B_{2}=C_{2}\cap {\rm Ker}(N_{2}^{3})={\rm Ker}(N_{2}^{2}).
$$
and if $C_{2}={\rm Ker}(N_{2}^{3})$, $B_{2}$ may vary freely in the 
projective line
$$
{\rm Ker}(N_{2})\subset B_{2}\subset {\rm Ker}(N_{2}^{3}).
$$
Therefore the $k'$-scheme of pairs $(B_{2}\subset C_{2})$ as 
above may be obtained by gluing the two projective lines over $k'$
$$
{\rm Ker}(N_{2})\subset B_{2}\subset {\rm Ker}(N_{2}^{3})
$$
and
$$
{\rm Ker}(N_{2}^{2})\subset C_{2}\subset {\rm Ker}(N_{2}^{4})
$$
along their marked points $B_{2}={\rm Ker}(N_{2}^{2})$ and $C_{2}={\rm 
Ker}(N_{2}^{3})$ respectively, and the stratum 
$(b_{1}^{+}=0,b_{2}^{+}=-1)$ is a ${\Bbb G}_{{\rm m}}$-torsor on this 
$k'$-scheme.

The stratum $(b_{1}^{+}=-1,b_{2}^{+}=0)$ is isomorphic to the stratum 
$(b_{1}^{+}=0,b_{2}^{+}=-1)$.

The stratum $(b_{1}^{+}=0,b_{2}^{+}=0)$ is equal to the $k'$-scheme of 
vector subspaces $B_{2}=C_{2}$ of dimension $3$ of $V_{2}$ such that
$$
{\rm Ker}(N_{2}^{2})\subset B_{2}=C_{2}\subset {\rm Ker}(N_{2}^{4})
$$
and is thus a projective line.

Similarly the stratum $(b_{1}^{+}=1,b_{2}^{+}=-1)$ and 
$(b_{1}^{+}=-1,b_{2}^{+}=1)$ are projective lines over $k'$.
\vskip 1mm

We may summarize this discussion by saying that the set $X(k')$ has
$$
q'(q'-1)q'+(q'-1)(2q'+1)+(q'-1)(2q'+1)+(q'+1)+(q'+1)+(q'+1)=q'^{3}+3q'^{2}+q'+1
$$
elements and that the set $X'(k')$ has
$$
(q'-1)(2q'+1)+(q'+1)+(q'+1)=2q'^{2}+q'+1
$$
elements, so that
$$
|X(k')|-|X'(k')|=q'^{2}(q'+1).
$$
\vskip 5mm
\centerline{\bf References}
\vskip 5mm
\parindent=10mm
\newtoks\ref
\newtoks\auteur
\newtoks\titre
\newtoks\editeur
\newtoks\annee
\newtoks\revue
\newtoks\tome
\newtoks\pages
\newtoks\reste
\newtoks\autre

\def\livre{\leavevmode\llap{[\the\ref]\enspace}%
\the\auteur\pointir
{\sl\the\titre},
\the\editeur,
{\the\annee}.
\smallskip\filbreak}

\def\article{\leavevmode\llap{[\the\ref]\enspace}%
\the\auteur\pointir
\the\titre,
{\sl\the\revue}
{\bf\the\tome},
({\the\annee}),
\the\pages.
\smallskip\filbreak}

\def\autre{\leavevmode\llap{[\the\ref]\enspace}%
\the\auteur\pointir
\the\reste.
\smallskip\filbreak}

\ref={Gr}
\auteur={A. {\pc GROTHENDIECK}}
\reste={Formule des traces de Lefschetz et
rationalit\'e des fonctions L, S\'eminaire Bourbaki 1964/65, {\it in Dix
expos\'es sur la cohomologie des sch\'emas}, North Holland, (1968), 31-45}
\autre

\ref={Ha}
\auteur={R. {\pc HARTSHORNE}}
\titre={Residues and Duality}
\editeur={Lecture Notes in Math. 20, Springer-Verlag}
\annee={1966}
\livre

\ref={Ho-Sh}
\auteur={R. {\pc HOTTA} and N. {\pc SHIMOMURA}}
\titre={The Fixed Point Subvarieties of Unipotent Transformations on 
Generalized Flag Varieties and the Green Functions}
\revue={Math. Ann.}
\tome={241}
\annee={1979}
\pages={193-208}
\article

\ref={Ka-Lu}
\auteur={D. {\pc KAZHDAN} and G. {\pc LUSZTIG}}
\titre={Fixed point varieties on affine flag manifolds}
\revue={Israel J. Math.}
\tome={62}
\annee={1988}
\pages={129-168}
\article

\ref={Ko}
\auteur={R. {\pc KOTTWITZ}}
\reste={Calculation of some orbital integrals, {\it in} R.P. Langlands and D.
Ramakrishnan (ed.), {\it The zeta functions of Picard modular surfaces},
Publ. CRM. Montreal, (1992), 349-362}
\autre

\ref={La-La}
\auteur={J.-P. {\pc LABESSE} and R.P. {\pc LANGLANDS}}
\titre={$L$-indistinguishability for ${\rm SL}(2)$}
\revue={Can. J. Math.}
\tome={31}
\annee={1979}
\pages={726-785}
\article

\ref={La-Sh}
\auteur={R.P. {\pc LANGLANDS} and D. {\pc SHELSTAD}}
\titre={On the definition of transfer factors}
\revue={Math. Ann.}
\tome={278}
\annee={1987}
\pages={219-271}
\article

\ref={Se}
\auteur={J.-P. {\pc SERRE}}
\titre={Corps locaux}
\editeur={Hermann}
\annee={1968}
\livre

\ref={SGA4}
\auteur={M. {\pc ARTIN}, A. {\pc GROTHENDIECK} and J.-L. {\pc VERDIER}}
\titre={Th\'eorie des Topos et Cohomologie Etale des Sch\'emas}
\editeur={Lecture Notes in Math. 269, 270, 305, Springer-Verlag}
\annee={1972/73}
\livre

\ref={SGA4${1\over 2}$}
\auteur={M. {\pc DELIGNE}, J.-F. {\pc BOUTOT}, A. {\pc GROTHENDIECK}, L. {\pc
ILLUSIE} and J.-L. {\pc VERDIER}}
\titre={Cohomologie Etale}
\editeur={Lecture Notes in Math. 569, Springer-Verlag}
\annee={1977}
\livre

\ref={Wa}
\auteur={van der {\pc WAERDEN}}
\titre={Algebra I}
\editeur={Springer-Verlag}
\annee={1971}
\livre

\vskip 10mm
\let\+=\tabalign

\line{\hbox{\kern 5mm{\vtop{\+ G\'erard LAUMON\cr
\+ URA 752 du CNRS\cr
\+ Universit\'e de Paris-Sud\cr
\+ Math\'ematiques, B\^at. 425\cr
\+ F-91405 ORSAY Cedex (France)\cr}}}
\hfill\hbox{{\vtop{\+ Michael RAPOPORT\cr
\+ Mathematisches Institut\cr
\+ Universit\"at K\"oln\cr
\+ Weyertal 86-90\cr
\+ D-50931 K\"OLN (Germany)\cr}}\kern 5mm}}

\bye